\documentclass{amsart}
\usepackage{graphicx,color}
\usepackage{url}
\usepackage{amsmath,amsfonts,amsthm}
\usepackage{amsmath}
\usepackage[psamsfonts]{amssymb}
\usepackage{amscd}
\usepackage{floatflt}
\usepackage[all,cmtip]{xy}
\newtheorem{thm}{Theorem}[subsection]

\theoremstyle{definition}

\newtheorem{con}[thm]{Conjecture}
\newtheorem{qst}[thm]{Question}
\newtheorem{problem}[thm]{Problem}

\theoremstyle{remark}

\numberwithin{equation}{section}

\newcommand{\Real}{\mathbb R}

\newcommand\bib[1]{\bibitem[#1]{#1}}
\newcommand{\weg}[1]{}

\def\margin#1{ \vspace{.5cm}  \hspace{-2cm}�\parbox{14cm}{  \texttt{ #1 }}   \vspace{.5cm}}%
\begin{document}

\title{Open problems and questions about geodesics }%
\author{Keith Burns}  %
\address{Department of Mathematics, Northwestern University Evanston, IL 60208-2730, USA}
\email{burns@math.northwestern.edu}
\author{Vladimir S. Matveev}
\address{Institut f\"ur Mathematik, Friedrich-Schiller-Universit\"at Jena, 07743 Jena,
  Germany}
\email{vladimir.matveev@uni-jena.de}


\begin{abstract}
The paper surveys open problems and questions related to geodesics defined by Riemannian, Finsler, semi Riemannian and magnetic structures on manifolds.
 It is an extended report on problem sessions held during the International Workshop on Geodesics in August 2010 at the Chern Institute of Mathematics in Tianjin.  

\end{abstract}
\maketitle

This paper is an extended report on problem sessions held during the International Workshop on Geodesics in August 2010 at the Chern Institute of Mathematics in Tianjin.  The focus of the conference was on geodesics in smooth manifolds. It was organized by Victor Bangert and Yiming Long and supported by AIM, CIM, and NSF.

\section{Notation and Definitions}
In this paper
$M$ is a connected $C^\infty$ manifold and $\widetilde M$ is its universal cover. We consider various structures on $M$ that create geodesics:
Riemannian (later also semi-Riemannian) metric $g$, Finsler metric $F$, magnetic structure~$\omega$.
 The unit tangent bundle defined by such a structure is denoted by 
 $SM$ and $\phi_t$ is the geodesic flow on $SM$.  The geodesic defined by a vector $v \in SM$ is denoted by $\gamma_v$; it is parametrized by arclength unless otherwise specified.

The free loop space $\Lambda M$ is the set of all piecewise smooth  mappings from the circle $S^1$ to the manifold equipped with the natural topology.  When $M$ has a Riemannian metric (or other suitable structure) we can consider the subspace $\Lambda_T M$ of loops with length $\leq T$. Similarly given two points $p,q \in M$, we denote by $\Omega(p,q)$ the space of piecewise smooth paths from $p$ to $q$ and by $\Omega_T(p,q)$ the subspace of paths with length at most $T$.

We define the energy functional on $\Lambda M$ by
$$
 c \mapsto \int_{S^1}\|\dot c(s)\|^2 \,d\lambda(s),
 $$ 
where $\|\cdot\|$ is the norm induced by a Riemannian metric or other suitable structure and $\lambda$ is Lebesgue measure on $S^1$ normalized to be a probability measure.
The critical points of this functional  are the closed geodesics for the metric (parametrized at the constant speed $Length(c))$. 
In the case of $\Omega(p,q)$, the circle  $S^1$ is replaced by the interval $[0,1]$.

\section{Closed geodesics} \label{mm}

Given a  Riemannian manifold, does there exist a closed geodesic? If yes, how many (geometrically different) closed geodesics must exist?  Existence is known if the manifold is closed.  For surfaces (and certain manifolds of higher dimension, e.g.\ $S^n$),  the result is essentially due to  Birkhoff ~\cite[Chap.~V]{Birkhoff1927}. 
 He used two arguments.  The first 
 is a variational argument in the free loop space $\Lambda M$.  
 From the subspace of $\Lambda M$ consisting of all loops homotopic  to a given loop, one selects a sequence of loops whose lengths converge to the infimum of the length function on the subspace.
   If the loops from this  sequence  lie in a compact region  of the manifold (for example, if the manifold itself is compact),  then by the Arzel\`a-Ascoli theorem 
  there exists a loop where the infimum is achieved; if the infimum is not zero, 
   this loop  is automatically a non-trivial closed geodesic. 
 These arguments show the existence of closed geodesics  in many cases, in particular for compact manifolds with non-trivial fundamental group.

If  the infimum of the length of the loops homotopic  to a given loop is zero, Birkhoff suggested another procedure (actually, a trick)  to prove the existence of closed geodesics, 
which is explained for example in   \cite[\S\S 6,7 of Chap. V]{Birkhoff1927}. This  trick 
 is nowadays called the \emph{Birkhoff minimax procedure}, and was used for example to show the existence of closed geodesics on any Riemannian sphere with dimension $\ge 2$.

 Let us recall this  procedure in the simplest case, when the manifold is $S^2$. We consider the foliation of the sphere
 without two points (north and south poles) 
  as in the picture:   if we think of the standard embedding in $\Real^3$, the fibers are the intersections of the sphere with the planes $\{x_3=\mathrm{const}\}$.

\begin{minipage}{\textwidth} {\includegraphics[width=.7\textwidth]{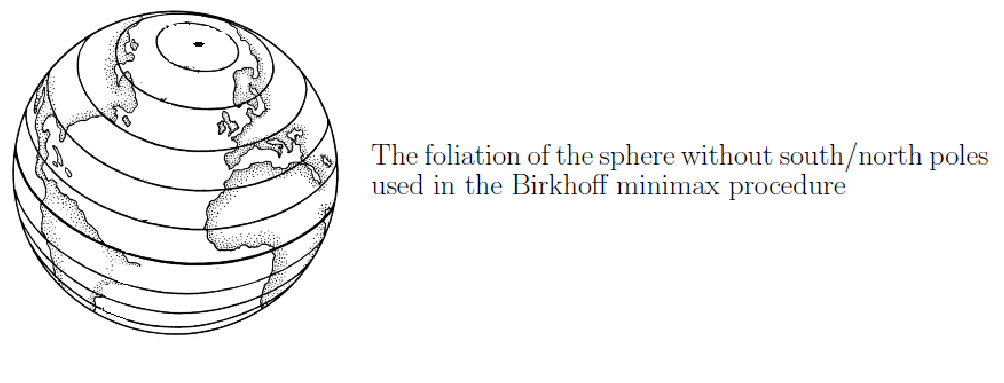}}

\end{minipage}

Now apply a curve shortening procedure  to every curve of this foliation, there are many options for such a  procedure. 
We need that the evolution of the curves in this procedure depends continuously on the curve. 
We obtain a sequence $F_i$ 
of the foliations of the sphere without two points into curves. For each foliation $F_i$  of the sequence,  let $\gamma_i$ be a leaf of maximal length. 
Because of topology, the lengths of the circles 
 $\gamma_i$   are bounded from below by a certain positive number.  By the Arzel\`a-Ascoli theorem, the sequence of curves $\gamma_i$ has a convergent subsequence. The limit of this subsequence is a 
 stable point of the shortening procedure, and is therefore a geodesic.

For arbitrary closed manifolds, the  existence of a closed geodesic is due to Lyusternik and Fet (\cite{Lyu-Fet1951}
and \cite{Fet1952}).   They  considered  the energy functional on the loop space $\Lambda M$ and   showed that the topology of $\Lambda M$ is complicated enough so that the energy functional must have critical points with non-zero energy (which are non-trivial closed geodesics).   Lyusternik and Schnirelmann \cite{Lyu-Sch1929a} and \cite{Lyu-Sch1929b} had earlier used a related argument to show that any Riemannian metric on $S^2$ has at least $3$ simple closed geodesics.

How many closed geodesics must exist on a closed manifold? For surfaces, the answer is known:  there are always infinitely many geometrically different closed geodesics. This is easily proved using Birkhoff's first argument when the fundamantal group is infinite. The remaining cases of the sphere and the projective plane were settled by \cite{Bangert1993} and  \cite{Franks1992}.  In higher dimensions  Rademacher has shown that a closed manifold with a generic Riemannian metric admits infinitely many geometrically different closed geodesics \cite{Rademacher1989}.

The main approach to proving the existence of infinitely many closed geodesics has been to apply Morse theory to the energy functional on the free loop space $\Lambda M$. The critical points of this functional are precisely the closed geodesics. 
But it has to be remembered that each  closed geodesic can be traversed an arbitrary number of times. It is thus important to distinguish geometrically different geodesics from repetitions of the same geodesic. This distinction is difficult to make and so far, despite some published claims, the existence of infinitely many closed geodesics on a general compact Riemannian manifold 
has not been proven. Nor have the resources of the Morse theory approach been fully exhausted.

One can modify these  classical questions in different directions; this will be done below. 

\subsection{Riemannian metrics on spheres }
Except in dimension 2, all that  is currently known for a general Riemannian metric on a sphere are the results that hold for all compact Riemannian manifolds.
\emph{Is there an extension of Lyusternik and Schnirelmann's result to higher dimensions, in particular to $S^3$}?

\subsection{Finsler metrics on $S^2$.}  \label{exkatok} 

The arguments of  Fet  (and Morse) can be adapted to the Finsler setting: 
one can show the existence of at least one closed geodesic on every compact manifold.  Moreover, if the Finsler 
 metric is reversible, or if the manifold is a surface other than the sphere or the projective plane, one can show the existence of infinitely many geometrically different closed geodesics. 
  
The following example constructed in \cite{Katok1973} shows that the number of closed geodesics on a 2-sphere  
with an irreversible metric can be~$2$.

Consider the sphere $S^2$ with the standard metric $g_0$ of constant curvature $1$  and the Killing field $V$ induced by the one parameter group of rotations about the north and south poles. Define a one parameter family of functions on $T^*M$ by
 $$
 F_\alpha(\xi) = \sqrt{g^*_0(\xi,\xi)} + \alpha \xi(V),
 $$
 where $g^*_0$ is the dual metric  of $g_0$. For small enough $\alpha \in \mathbb R$ there is a Randers Finsler metric, $F_\alpha$ such that
 $$
 \|v\|_\alpha = F_\alpha \circ L^{-1}_\alpha(v),
 $$
 where $L_\alpha: T^*M \to T^M$ is the Legendre transformation induced by the function $\frac12 F^2_\alpha$. If $V$ is scaled so that its integral curves have the same period as the geodesics of $g_0$, then $g_\alpha$ has precisely two closed geodesics for each irrational $\alpha$. They are the equator (the unique great circle tangent to $V$) traversed in both directions.
See \cite{Ziller1982} and \cite{Rademacher2004} for more details.


This  example can be generalized to higher dimensions:    one obtains a Finsler metric on $S^n$ with precisely $2[(n+1)/2]$ distinct prime closed geodesics; see \cite{Katok1973}.

 V.~Bangert and Y.~Long in \cite{Ban-Lon2010} and \cite{Long2006} proved that there are always at least 2 distinct prime closed geodesics for every irreversible 
    Finsler metric on $S^2$.

\begin{con}[Long, Bangert,   Problem 15 from \cite{Alvarez2006}] {\it Every irreversible Finsler metric on $S^2$ has either
exactly $2$ or infinitely many distinct prime closed geodesics.} \end{con}

 There exist  results supporting this conjecture. In particular,  
H.~Hofer, K.~Wysocki and E.~Zehnder in \cite{Ho-Wy-Ze2003} studied Reeb orbits on contact $S^3$.
Their result can be projected down to $S^2$, and implies that the total number of distinct prime closed geodesics for a bumpy Finsler metric on $S^2$ is either 2
or infinite, provided the stable and unstable manifolds of every hyperbolic closed geodesics
intersect transversally. See also  
  \cite{Har-Pat2008,Rademacher2016}.

The closed geodesics in Katok's example are elliptic.  

\begin{con}[Long]  {\it The existence of one hyperbolic prime closed geodesic on a Finsler
$S^2$ implies the existence of infinitely 
  many distinct prime closed geodesics.} \end{con} 

\begin{con}[Long]  {\it  Every Finsler $S^2$ has at least one elliptic prime closed geodesic.} \end{con}

The conjecture   agrees with  a result of  
Y.~Long and W.~Wang who   proved that there are always at least 2 elliptic prime
closed geodesics on every irreversible       Finsler $S^2$,  if the total number of prime closed geodesics is finite~ \cite{Lon-Wan2008}. 
The conjecture does not contradict  \cite{Grjuntal1979} where an example of a metric such that all  closed simple geodesics are   hyperbolic is constructed.  Indeed, for certain metrics on the sphere 
 (and even for the metric of certain ellipsoids) most prime closed geodesics are not simple.

\subsection{ Of complete Riemannian metrics  with finite volume}
\

\begin{qst} [Bangert]
{\it Does every complete Riemannian manifold with finite volume have at least one closed geodesic?}\end{qst}

The question  was   answered affirmatively for dimension 2. Moreover, in dimension 2 a complete Riemannian manifold of finite volume even has infinitely many geometrically different geodesics  \cite{Bangert1980}. The argument that was used in the  proof is based on the Birkhoff minimax procedure 
we recalled in the beginning of \S \ref{mm}, and does not work in dimensions $\ge 3$. One  can even hope to construct counterexamples in the class of Liouville-integrable geodesic flows.
In this case,  most orbits of the geodesic flow are rational or irrational windings on the Liouville tori; they are closed, if all of the  corresponding frequencies are  rational. 
Since there are essentially  $(n-1)$ frequencies  in dimension $n$, one can hope that if $n > 2$ it would be  possible to ensure that there is always at least one irrational frequency.  
Initial  attempts to find    a counterexample on  $T^2\times \Real$ with the metric of the form $a(r)d\phi^2 + b(r) d\psi^2+ dr^2 $  were, however, unsuccessful.

A more difficult problem would be to prove {\it that on any  complete Riemannian manifold of finite volume
there exist  infinitely many closed geodesics}.

Questions of this nature are also interesting in  the realm of Finsler geometry. One would expect the results for reversible Finsler metrics to be very similar to those for Riemannian metrics. For irreversible Finsler metrics, Katok's example shows that there are compact Finsler manifolds with  only finitely many closed geodesics, but it is still possible that all non-compact Finsler manifolds with finite volume might have infinitely many closed geodesics. On the other hand, the following question is completely open:

\begin{qst}[Bangert]  {\it Does there exist an  irreversible Finsler metric  of finite volume
on $\Real\times S^1$ with no  closed geodesics?} \end{qst}


\subsection{Of  magnetic flows  on closed surfaces.}

It is known that the trajectory of a charged particle  in the presence of magnetic forces  (=``magnetic geodesic'') is described by a Hamiltonian system with the ``kinetic"  Hamiltonian of the form $\sum_{i,j} p_ip_jg^{ij}$ on $T^*M$ with the symplectic form $dp\wedge dx +\pi^* \omega$, where $\omega $ is a  closed (but not necessarily exact) form on $M$ and $\pi:T^*M\to M$ is the canonical projection.

We assume that our surface $M^2 $ is closed.

\begin{qst}[Paternain]   {\it Is there at least one closed magnetic geodesic in every positive energy level?} \end{qst}

The only known case in which the answer is negative is when $\omega$ is the area form for a hyperbolic surface and the magnetic flow is the horocycle flow. On the other hand Ginzburg showed that there is always a closed orbit if the magnetic field is weak enough
\cite{Ginzburg1989, Ginzburg1996}.

If $\omega$ is exact, the affirmative answer for all closed orientable surfaces and all energy levels was obtained   in 
\cite{Co-Ma-Pa2004}, which is partially based on \cite{Taimanov1992}.    Actually, by \cite{Ab-Ma-Pa2015,AMMP2017},  in this case there are at  least two  closed geodesics on almost every energy level. We give some details.

If the energy is greater than Ma\~n\'e's critical value $c_u$, then the magnetic flow
is conjugate to a Finsler geodesic flow, which has $\geq 2$ closed orbits in the case of $S^2$ (by \cite{Ban-Lon2010}) and infinitely many closed orbits if the genus of the surface is positive (at least two for each free homotopy class). The lower bound of 2 is sharp because of the Katok example mentioned above.

For energy less than $c_u$, the results about more than one closed geodesic are only for almost every energy level.
Contreras obtained at least 2 closed geodesics \cite{Contreras2006}. In \cite{Ab-Ma-Pa2015} this was improved to at least 3, and in fact infinitely many under 
a non-degeneracy assumption, which was removed in \cite{AMMP2017}.

In the case of magnetic monopoles (i.e.\ when $\omega$ is not exact), the answer is affirmative if the energy level is contact; this is an immediate corollary of the proof of the Weinstein conjecture in dimension 3 by Taubes \cite{Taubes2007}. On the other hand, Benedetti has given examples of magnetic monopoles on $S^2$ with energy levels that are not contact; see Proposition 4.9 in \cite{Benedetti2016}.  On $S^2$, if the form $\omega$ is nowhere vanishing, the existence of a closed magnetic geodesic on sufficiently big energy level is established in \cite[Corollary 4.6]{Ass-Ben2016}; see also \cite{Ginzburg1989,Ginzburg1996}.

\subsection{Kropina metrics with  nonintegrable distribution}

A Kropina metric on a Riemannian manifold $(M,g)$ is the function on the tangent bundle given by $F(\xi)= \tfrac{g(\xi, \xi)}{\omega(\xi)}$, where 
$\omega$ is a nowhere vanishing 1-form. We will assume that $\omega\wedge d\omega$ nowhere vanishes;  this   assumption implies  that the distribution $\textrm{Ker}(\omega)$ is not integrable.  A  Kropina metric is not a Finsler metric, since it is not defined on the kernel of $\omega$, but still the function $F$ is one-homogeneous and strictly convex at $\xi\not\in \textrm{Ker}(\omega)$; one can use it to define the distance on $M$ and geodesics on $M$.

Recent interest in Kropina metrics is due to their relations to the so-called chains in CR geometry, see \cite{CMMM2019}.

 \begin{qst} \label{kr1}   {\it Does every Kropina metric on a closed manifold  have a   closed geodesic?} \end{qst} 
 
 If  the fundamental group is nontrivial, the answer is positive  and the proof is essentially the same as in  the Riemannian case. But the Birkhoff-Morse-Fett approaches that  were used to prove the existence of closed  geodesics for simply-connected Riemannian manifolds can not be generalized directly, in particular because of the existence of examples of Kropina metrics whose geodesics of length $1$  lie in an  arbitrarily small neighborhood.

As a simple   special case   of Question \ref{kr1}  we suggest  considering  $M= S^3= \textrm{Spin}(3)$ equipped 
 with arbitrary $g$ and  a 1-form $\omega$ that is
homogeneous with respect to the left action of $\textrm{Spin}(3)$.

\section{Path and loop spaces} \label{path}

\weg{
Let us recall the ``Morse approach'' to proving the existence of closed geodesics on closed manifolds, and to estimating the number of closed geodesics of length less than certain number. 
Let $(M, g)$ be a Riemannian manifold.  
Recall that the {\it loop space} is the space of all (smooth, continuous, or piece-wise smooth) mappings $\gamma:S^1\to M$ equipped with the natural topology. On the loop space, let us consider the energy functional: 

$$\gamma\mapsto \int_{S^1}g( \dot \gamma, \dot\gamma).$$ 

The critical points of this functional  are the closed geodesics for the metric. 
If the topology of the loop space is sufficiently complicated, every function on this space must have  critical points 
lying on  different levels.  Thus, studying the topology of the loop space may allow one 
 to show the existence of closed geodesics. 
 
 One can also consider the space of the loops of length $\le T$ for a certain number $T$. 
  Its topology is then useful to 
   estimate the number of closed geometrically different geodesic loops of length $\le T$. 
  
One can also replace the space of loops by the space of paths of length $\le T$   connecting two points;
  the topology of the path spaces is     related to the number of  geodesics of length less than a certain number $T$ connecting two points. 
}

As noted in the previous section, one of the main approaches to proving the existence of closed geodesics is to use topological complexity of the loop space $\Lambda M$ to force the existence of critical points of the energy functional. Loops with length $\leq T$ correspond to critical points  in $\Lambda_T M$. Similarly geodesics joining two points $p$ and $q$ can be studied by investigating the path spaces $\Omega(p,q)$ or $\Omega_T(p,q)$. The homology of these spaces have been much studied.

\subsection{ Sums of the Betti numbers }

Let $p, q$ be points in a Riemannian manifold. The space $\Omega_T(p,q)$ of paths from $p$ to $q$ with length $\leq T$ has the homotopy type of a finite complex (see e.g.\  \cite{Milnor1963}), 
 and hence the sum of its Betti numbers is finite for each $T$. The same is true for the space $\Lambda_TM$ of loops with length at most~$T$.

 \begin{qst}[Paternain]    {\it How do the sums of the Betti numbers for $\Lambda_TM$ and $\Omega(_T(p,q)$  grow as $T\to \infty$? Does the growth  depend on the metric?} \end{qst}

  It was shown by Gromov in \cite{Gromov2001} that   
 if $M$ is simply connected
there is a constant $C$ such that the
sum of the Betti numbers of $\Lambda_{CN}$  is at least
$
\sum_{i=0}^N b_{i}(\Lambda)$.

 \begin{qst}[Gromov] {\it 
Does the number of closed geodesics with length $\leq T$ grow exponentially as $T \to \infty$ if the Betti numbers of the loop space grow exponentially?}\end{qst}

From \cite{Gromov1978} it follows that the answer is positive  for generic metrics.

Here is a potentially interesting example. Consider $f:S^3\times S^3\to S^3\times S^3$ such that  the induced action on 
$H_3(S^3\times S^3)$  is hyperbolic. Let $M$ be the mapping torus for $f$, i.e.~ $S^3\times S^3 \times [0,1]$ with $(x,1)$ identified with $(f(x),0)$ for each $x \in S^3 \times S^3$. Then $\pi_1({M})=Z$  and the universal cover $\widetilde  M$ is homotopy equivalent to $S^3\times S^3$. Hence the Betti numbers of the loop space grow polynomially. On the other hand, the hyperbolic action of $f$ on $H_3$ gives hope for exponential growth of the sum of the Betti  numbers of $\Omega_T(p,q)$.

\subsection{Stability of minimax   levels (communicated by Nancy Hingston).}

Let $M$ be a compact manifold with a  Riemannian (or Finsler) metric $g$. \
Let  $G$ be a finitely generated
abelian group. \ Given a non-trivial homology class 
$X \in H_{\ast}(\Lambda M;G)$, the critical (minimax) level of $X$ is
\begin{align*}
cr\,X &  =\inf\{a:X\in\text{Image }H_{\ast}(\Lambda_a M;G)\}\\
&  =\underset{x\epsilon X}{\inf}\text{ }\underset{\gamma\in\text{Image}%
x}{\sup}Length({\gamma})
\end{align*}
Here $\Lambda_a M$ is the subset of the free loop space $\Lambda M$
consisting of loops whose  length is 
at most $a$. 
The second definition is the minimax definition: the
singular chain $x$ ranges over all representatives of the homology class $X$,
and  $\gamma$ over all the points in the image of $x$, which are loops in $M$.

\bigskip

\begin{qst}[Hingston] {\it Do there exist a metric on }$S^{n}$\textit{ and a
homology class }$X \in H_{\ast}(\Lambda M;%
\mathbb{Z}
)$\textit{ with }%
\[
0 < cr(mX) < c r(X)
\]
\textit{for some }$m \in
\mathbb{N}
$?  \end{qst}

The special case of this question, when dimension is two and $m$ is odd,  was answered negatively in \cite[Corollary 1.3]{Cha-Lio2014}.

\bigskip

Given a metric $g$ on $M$ and a finitely generated abelian group $G$, the
\textit{global mean frequency} is defined as
\begin{equation} \label{ast}
\alpha_{g,G}=\underset{\deg X\rightarrow\infty}{\lim}\frac{\deg X}{crX},
\end{equation}
where the limit is taken over all non-trivial homology classes $X \in
H_{\ast}(\Lambda S^{n};G)$.

\bigskip The Resonance Theorem from \cite{Hin-Rad2013}  says that if $M$ is a sphere and $G$ is a
field, the limit \eqref{ast} exists. It is clear in this case that
$\alpha_{g,G}$ depends on the metric $g$.  But does it really depend on the
field $G$? \ \ The degree of $X$ does not depend on anything but $X$. \ But
what about the critical level $cr\,X$? \ Does $cr\,X$ depend on the coefficients? \

Let us note that 
for the spheres the non-trivial homology groups of the free loop space (with
integer coefficients) are all $\mathbb{Z}$ or $%
\mathbb{Z}
_{2}=:%
\mathbb{Z}
/2%
\mathbb{Z}
$. \ (For odd spheres they are all $%
\mathbb{Z}
$.) \ \ Let us look at the case where $X \in H_{k}(\Lambda;%
\mathbb{Z}
)=%
\mathbb{Z}
$. \ \ \ For each $m \in
\mathbb{N}
$ there is a critical level
\[
cr(mX)=\inf\{a:mX \in \text{Image }H_{\ast}(\Lambda_{ a}M)\}.
\]
If $j,m \in
\mathbb{N}
$, then clearly (since Image $H_{\ast}(\Lambda_{a})$ is an additive
subgroup of $H_{\ast}(\Lambda)$)%

\[
cr(jmX)\leq cr(mX)\text{.}%
\]
But can there be strict inequality?  Here is a little
``example'' to show how this would affect the global mean frequency : Suppose it
were the case that there were real numbers $a<b<c<d$ with%

\[
cr(mX)=\left\{
\begin{array}
[c]{c}%
d\text{ if }\gcd(m,3)=\gcd(m,7)=1\\
c\text{ if }3|m\text{ \ but }\gcd(m,7)=1\\
b\text{ if }7|m\text{ \ but }\gcd(m,3)=1\\
a\text{ if }21|m
\end{array}
\right\}
\]

\begin{con}[Hingston]
\begin{align*}
\alpha_{g,G}  &  =a\text{ if }G=\mathbb{Q}\\
\alpha_{g,G}  &  =b\text{ if }G = \mathbb{Z}_{3}\\
\alpha_{g,G}  &  =c\text{ if }G = \mathbb{Z}_{7}\\
\alpha_{g,G}  &  =a\text{ if }G= \mathbb{Z}_{p},p\neq3,7.
\end{align*} \end{con}

\section{Curvature conditions and hyperbolicity of the geodesic flow}
Many results about geodesics and the geodesic flow assume that all sectional curvatures are negative or one of the following increasingly weaker properties:
\begin{enumerate}
\item all sectional curvatures are non-positive,
\item no focal points,
\item no conjugate points.
\end{enumerate}
These properties can be characterized by the behaviour of Jacobi fields:
\begin{enumerate}
\item[(0)] Negative curvature: the length of any (non-trivial) Jacobi field orthogonal to a geodesic is a strictly convex function.

\item Non-positive curvature: the length of any Jacobi field is a convex function.

\item No focal points: the length of an initially vanishing Jacobi field is a 
 non-decreasing function along a geodesic ray.

\item No conjugate points: a (non-trivial) Jacobi field can vanish at most once.
\end{enumerate}

The no focal point property is equivalent to convexity of spheres in the universal cover. Most interesting results about manifolds with non-positive curvature extend readily to manifolds with no focal points.

For a compact manifold, negative curvature implies that the geodesic flow is uniformly hyperbolic, in other words an Anosov flow. This means that there is a 
$D\phi_t$-invariant splitting of the tangent bundle of the unit tangent bundle $SM$,
 $$
 TSM = E^s \oplus E^0 \oplus E^u,
 $$
 in which $E^0$ is the one dimensional subbundle tangent to the orbits of the geodesic flow, and there are constants $C \geq 1$ and $\lambda > 0$ such that for any $t \geq 0$ and any vectors $\xi \in E^s$ and $\eta \in E^u$ we have
 \begin{equation} \label{hyp}
 \|D\phi_t(\xi)\| \leq Ce^{-\lambda t}\|\xi\| 
 \quad\text{and}\quad
 \|D\phi_{-t}(\eta)\| \leq Ce^{-\lambda t}\|\eta\|.
 \end{equation}
(Here we have in mind the usual Sasaki metric; the same property would also hold for any equivalent metric with different constants $C$ and $\lambda$.) This splitting is H\"older-continuous, but usually not smooth.

The bundles $E^s$ and $E^u$ for an Anosov geodesic flow are integrable; their integral foliations are usually denoted by $W^s$ and $W^u$. 
 The lifts to the universal cover of the leaves of $W^s$ and $W^u$ are closely related to horospheres. If $\widetilde v$ is the lift to $S\widetilde M$ of $v \in SM$, the lifts to $S\widetilde M$ of $W^s(v)$ and $W^u(v)$ are formed by the unit vectors that are normal to the appropriate horospheres orthogonal to  $\widetilde v$ and  are on the same side of the horosphere as $\widetilde v$; see the picture below.

{\includegraphics[width=.5\textwidth]{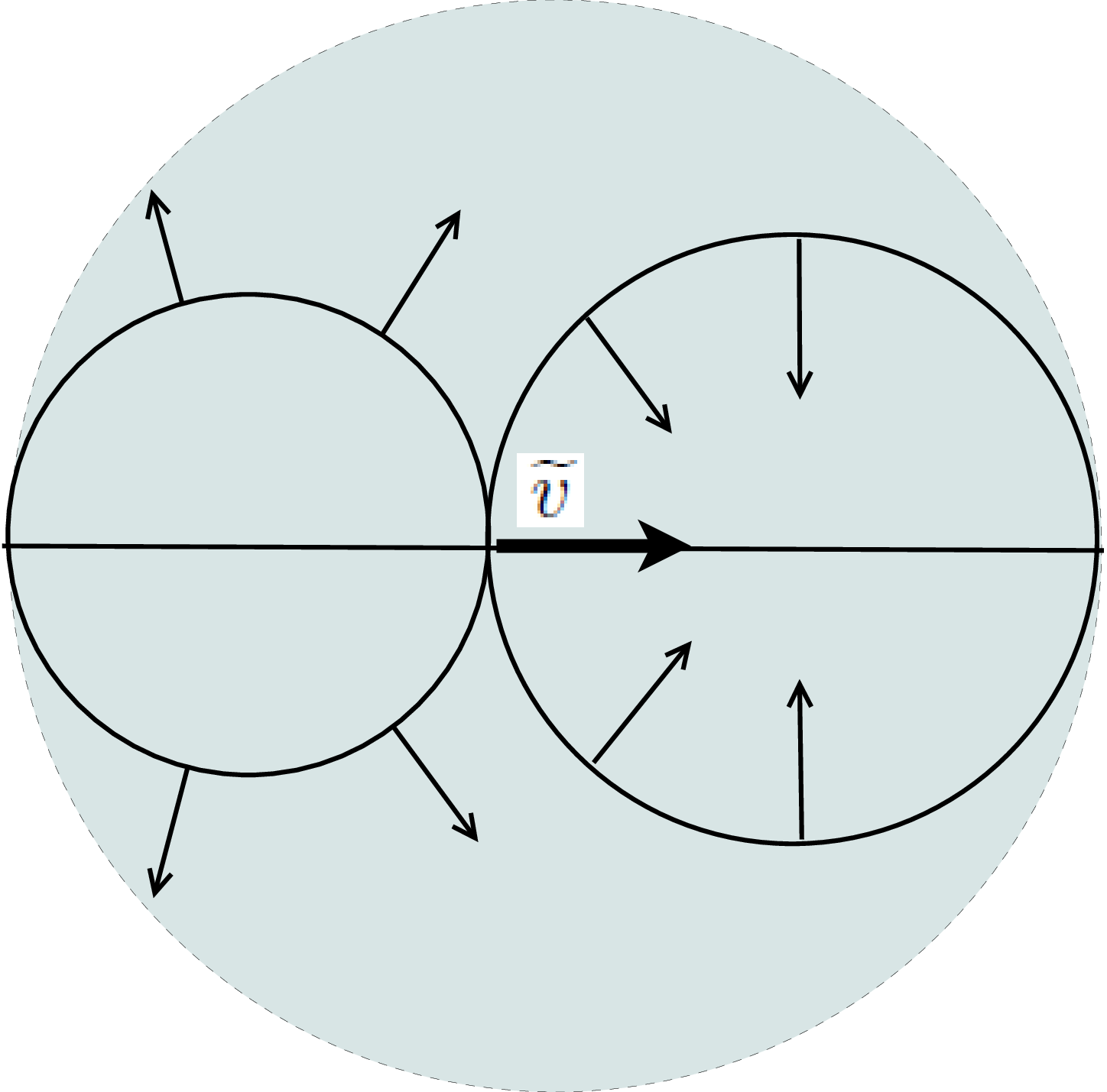}}

\subsection{Relations between these concepts} \label{relationsbetween}

The notions introduced above are related as follows:

\[
\begin{xy}
\xymatrixcolsep{1pc}
\xymatrixrowsep{2pc}
\xymatrix{
{\textrm{negative curvature}}\ar@2[r]\ar@2[d]&{\textrm{non-positive curvature}}\ar@2[r]& 
{\textrm{no focal points}}\ar@2[d]\\
{\textrm{Anosov geodesic flow}}\ar@2[rr]&&{\textrm{no conjugate points}} \\
}
\end{xy}
\]
That Anosov geodesic flow implies no conjugate points is proved in part B of \cite{Mane1987}.

The class of compact manifolds that support metrics with variable negative curvature  is much larger than the class that support hyperbolic metrics (of constant negative curvature). The earliest  examples of  manifolds that support variable but not constant negative curvature were given by Mostow-Siu \cite{Mos-Siu1980} and Gromov-Thurston \cite{Gro-Thu1987}. Recent work  of Ontaneda has vastly increased the supply of examples \cite{Ontaneda2011, Ontaneda2014}.

 \begin{con}[Klingenberg] {\it If a closed manifold   admits a metric with  Anosov geodesic flow, then it admits a metric with  negative  sectional curvature.}  \end{con}

All known examples of metrics with Anosov geodesic flow are perturbations of metrics with negative curvature.
It is not difficult to show using  the uniformization theorem and the Thurston geometerization  theorem that conjecture is true in dimensions 2 and 3. Klingenberg \cite{Klingenberg1974} showed that seven properties of Riemannian manifolds with negative curvature extend to those with Anosov geodesic flow. One of these is Preissman's theorem \cite{Preissman1943} that the fundamental group of a manifold with negative sectional curvatures cannot contain a copy of ${\mathbb Z} \times  {\mathbb Z}$.

Several examples of manifolds that admit metrics of non-positive curvature but cannot support a metric of negative curvature (or with Anosov geodesic flow) can be found in the introduction to  \cite{Ba-Br-Eb1985}.  These examples have a copy of ${\mathbb Z} \times  {\mathbb Z}$  in their fundamental group. The simplest of them is due to Heintze who took  two cusped hyperbolic 3-manifolds, cut off the cusps, glued the two pieces together along their boundary tori and then smoothed out the metric to obtain a manifold with 
 non-positive curvature.

\begin{qst}  {\it If $(M,g)$ has no conjugate  points, does $M$ carry a metric 
 with  non-positive sectional curvature?}  \end{qst} 

Again it is not difficult to use the uniformization   theorem  to show that the answer is affirmative in dimension 2.
But the problem is still open even in dimension 3. See \cite{Cro-Sch1986}, \cite{Lebedeva2002} and  \cite{Iva-Kap2014} for results of the nature that the fundamental groups of manifolds with no conjugate points share properties with those of non-positive curvature. Section 8 of \cite{Iva-Kap2014} contains a number of open problems of which we mention:

\begin{qst}  {\it Is the fundamental group of a closed manifold without conjugate points semihyperbolic?}  \end{qst} 

Semihyperbolicity is a condition introduced by Alonso and Bridson \cite{Alo-Bri1995} to describe non-positive curvature in the large for an 
arbitrary metric space.

The geodesic flow is partially hyperbolic if there are a
$D\phi_t$-invariant splitting
 $$
 TSM = E^s \oplus E^c \oplus E^u
 $$
 and constants $C \geq 1$ and $\lambda > \mu > 0$ such that (\ref{hyp}) holds and in addition for all $t$ and all $\zeta \in E^c$ we have
 $$
  C^{-1}e^{-\mu|t|}\|\zeta\| \leq \|D\phi_t(\zeta)\| 
      \leq Ce^{\mu|t|}\|\zeta\|.
 $$

 Anosov geodesic flows are degenerate examples of partially hyperbolic geodesic flows. Genuine examples (with $\dim E^c \geq 2$)  have been constructed by Carneiro and Pujals; see \cite{Car-Puj2011} and \cite{Car-Puj2013}. They deformed a higher rank locally symmetric space of non-compact type. This answered a question raised by Herman in emails to many dynamicists.

\section{Negative curvature and hyperbolicity of the geodesic flow}

\subsection{ Does the marked length spectrum  determine the metric?}

Suppose $(M, g)$ is a closed Riemannian manifold. The marked length spectrum of $M$ is the
function  
that assigns to each free homotopy class of loops the
infimum of the length of loops in the class (i.e.~the length of the  shortest closed geodesic lying in this class).

 \begin{con}[\cite{Bur-Kat1985}] \label{mlscon}
{\it Two metrics with negative curvature on a compact manifold must be isometric if they have the same marked length spectrum.} 
\end{con} 
Croke and Otal showed that this conjecture is true for metrics on surfaces \cite{Otal1990,Croke1990}. Indeed it is enough to assume that the metrics have
 non-positive curvature. The problem is open in higher dimensions. Hamenst\"adt showed that the geodesic flows of the two metrics must be $C^0$-conjugate, thereby reducing the problem to Conjecture~\ref{conjcon} in the next subsection \cite{Hamenstadt1991}. Guillarmou and Lefeuvre \cite{Gui-Lef2018} have shown that the conjecture is locally true in sense that two close enough metrics with Anosov geodesic flow and the same length spectrum are isometric.

One obtains  a natural interesting modification  of this conjecture by replacing negative curvature by non-existence of conjugate points.  On the other hand some restriction on the metrics is necessary as the examples of Croke and Kleiner which we describe in the next subsection provide non-isometric metrics with the same marked length spectrum. It is also necessary to consider the marked length spectrum rather than the length spectrum. 
Vign\'eras gave examples of non-isometric hyperbolic surfaces that have the same set of lengths for their closed geodesics \cite{Vigneras1980}.

One can ask a similar question for non-manifolds. 
Consider a 2-dimensional (metrical) simplicial complex such that  every simplex is hyperbolic with geodesic edges   and such that $CAT(-1)$ condition holds on every vertex. Assume in addition that every edge is contained in at least two simplices.

\begin{qst}[Benjamin Schmidt] {\it Does the marked length spectrum determine such a metric (in the class of all locally $CAT(-1)$ metrics on this space, or in the class of all 2-dimensional  (metrical) simplicial  complexes
with the above properties
 homeomorphic to the given complex)?} \end{qst}

If the simplicial complex is (topologically) a manifold, the answer is positive and is due to   
 \cite{Her-Pau1997}. One can also ask the question in higher dimensions.
An easier, but still interesting  version of the question  is when we assume that every edge is contained in at least three simplices.

The questions above  are  closely related to the boundary rigidity problem, which we now recall.  
Given a compact manifold  $M$  with smooth boundary $N$, a
Riemannian metric $g$ on $M$ induces a non-negative real valued function  $d$ 
 on
$N \times  N$ where $d(p, q)$ is the distance in $(M, g)$ between $p$ and $q$.  
We call $(M,g)$  boundary rigid if a Riemannian manifold $g'$ on $M$
that induces the same function on $N \times N$  must be isometric to $g'$.

\begin{qst}[Boundary rigidity problem] \label{brquest}
{\it What conditions on $M$, $N$ and  $g$ 
 imply   boundary rigidity?} 
\end{qst} 

Michel conjectured in \cite{Michel1981} that $(M,g)$ should be boundary rigid if it is \emph{simple}. Simplicity means that $M$ is simply connected and any two points of $M$ are connected by a unique geodesic segment $\gamma_{p,q}:[0,1] \to M$ and this unique geodesic has the property that $\gamma_{p,q}(t) \notin N$ if $0 < t < 1$. Michel’s conjecture is known to be true when $(M,g)$ is a simple subspace of a Euclidean space \cite{Gromov1983}, 
when $(M.g)$ is a simple subspace of a symmetric space with negative curvature 
\cite{Be-Co-Ga1995}, and when $M$ is two dimensional \cite{Pes-Uhl2005}.

One can   modify Question~\ref{brquest}  by requiring that the other Riemannian metric $g'$ in the definition of boundary rigid manifolds above also satisfies some additional assumption. 
Special cases of  the last question were answered in \cite{Croke1990,Cr-Da-Sh2000,Sha-Uhl2000,Bur-Iva2013}.

One can also ask the boundary rigidity question for Finsler metrics; recent references  with non-trivial results include 
\cite{Coo-Del2010, Bur-Iva2010}.

\subsection{The conjugacy problem}
Two Riemannian manifolds $(M_1,g_1)$ and $(M_2,g_2)$ have $C^k$-conjugate geodesic flows if there is an invertible map $h: S^M_1 \to S^M_2$ such that $h$ and $h^{-1}$ are $C^k$ and
$h \circ \phi^1_t = \phi^2_t \circ h$ for all $t$. Here $\phi^i_t$ is the geodesic flow for the metric $g_i$. A long standing conjecture is

\begin{con} \label{conjcon}
Compact Riemannian manifolds with negative curvature  must be isometric if they have $C^0$-conjugate geodesic flows.
\end{con}

As with Conjecture~\ref{mlscon} it may be possible to relax the hypothesis of negative curvature to no conjugate points, but some restriction on the metrics is required.
Weinstein pointed out that all of the 
  Zoll metrics on $S^2$ have conjugate geodesic flows. 
  Using the explicit rotationally-symmetric examples of Zoll metrics (Tannery metrics in the terminology of  \cite[Chapter 4]{Besse1978}) Croke and Kleiner constructed examples of different  metrics on an 
  arbitrary closed 
  manifold such that their geodesic flows are conjugate and  their marked length spectra coincide, see \cite{Cro-Klei1994}. 
  
  The conjecture is true in dimension two; indeed it is enough to assume that one of the manifolds has non-positive curvature \cite{Cr-Fa-Fe1992}. It follows from the theorem of Besson-Courtois-Gallot discussed in the next subsection that the conjecture holds if one of the metrics
  is locally symmetric.
  Croke and Kleiner showed that $C^1$ conjugacy implies that the volumes of the manifolds are the same with no restrictions on the metrics 
  and implies isometry of the metrics if one of them has a global parallel vector (as is the case in a Riemannian product) \cite{Cro-Klei1994}. 
  
\subsection{Entropy and locally symmetric spaces with negative curvature}

Let $h_{Liou}$ denote the entropy of the geodesic flow with respect to the Liouville measure, $h_{top}$ its topological entropy and $h_{vol}$ the volume entropy 
 $$
 h_{vol} = \lim_{r \to \infty} \frac1r \log \text{\rm Vol}\,B_{\widetilde M}(p,r),
 $$
 where $B_{\widetilde M}(p,r)$ is the ball of radius $r$ around a point $p$ in the universal cover $\widetilde M$. It is always true that $h_{vol} \leq h_{top}$ and they are equal if there are no conjugate points, in particular if the curvature is negative \cite{Manning1979, Fre-Man1982}.
 By the variational principle for entropy $h_{Liou} \leq h_{top}$ for any metric; these entropies are equal for a locally symmetric metric.
 
Katok \cite{Katok1982} showed that if $g$ is an arbitrary metric and $g_0$ a metric of constant negative curvature on a surface of genus $\geq 2$ such that $\text{area}(g) = \text{area}(g_0)$, then
 $$
 h_{Liou}(g) \leq h_{Liou}(g_0) =h_{top}(g_0) \leq h_{top}(g)
 $$
 and both inequalities are strict unless $g$ also has constant negative curvature. He explicitly formulated the following influential and still open conjecture:
 
 \begin{con}[\cite{Katok1982}]\label{Katokentropycon}
 $h_{top} = h_{Liou}$ for a metric of negative curvature on a compact manifold if and only if it is locally symmetric.
 \end{con}
 
Katok's arguments extend to higher dimensions provided the two metrics are conformally equivalent.
 Flaminio showed that Conjecture~\ref{Katokentropycon} holds for deformations of constant curvature metrics \cite{Flaminio1995}.

 Katok's paper implicitly raised the questions of whether  in higher dimensions $h_{top}$ is minimized and $h_{Liou}$ maximized (among metrics of negative curvature of fixed  volume on a given manifold) by the  locally symmetric metrics. It is now known that the locally symmetric metrics are critical points for both entropies \cite{Ka-Kn-We1991}, but   the locally symmetric spaces  do not maximize $h_{Liou}$ except in the two dimensional case \cite{Flaminio1995}. The topological entropy, however, is  minimized.
 
 Gromov \cite{Gromov1983} conjectured that if $f: (M,g) \to (M_0,g_0)$ is a continuous map of degree $d \neq 0$ between compact connected oriented $n$-dimensional Riemannian manifolds and $M_0$ has constant negative curvature, then
 $$
  h_{vol}(g)^n\text{Vol}(M,g) \geq |d|  h_{vol}(g_0)^n\text{Vol}(M_0,g_0)
  $$
  with equality if and only if $g$ also has constant negative curvature and $f$ is homotopic to a $d$-sheeted covering. The quantity $h^n\text{Vol}$ has the advantage of being invariant under homothetic rescaling of the manifold; this obviates the assumption that the two manifolds have the same volume.
  
  \begin{thm}[\cite{Be-Co-Ga1995}] \label{BCGthm}
  The above conjecture of Gromov is true even when $(M_0,g_0)$ is allowed to be a locally symmetric space of negative curvature.
  \end{thm}
  
  We refer the reader to Section 9 of \cite{Be-Co-Ga1995} and the survey \cite{Eberlein2001} for a list of the many corollaries of this remarkable theorem.
  
  \subsection{Regularity of the Anosov structure}
  An extensive list of results about the regularity of the Anosov splitting for the geodesic flow in negative curvature (and for other Anosov systems) can be found in the introduction to \cite{Hasselblatt1994}. In dimension two, or if the curvature is $1/4$-pinched, the splitting is always at least $C^1$. In higher dimensions it is always H\"older continuous but typically not $C^1$. It is natural to ask:

\begin{qst} \label{C^1quest}
{\it Does there exist a closed Riemannian manifold of negative curvature such that $E^s $ and $ E^u$ are $C^1$ smooth, but the dimension is $> 2$ and the metric is not $1/4$-pinched}? \end{qst} 

Another motivation for this question is the question of Pollicott in the next subsection.

  A program to construct open sets of metrics whose geodesic flows have Anosov splitting that have low regularity of the Anosov splitting on large subsets of the unit tangent bundle is explained in section 4 of \cite{Has-Wil1999}; see in particular Proposition 12. The idea is to perturb a suitable base example. This base example must have directions in its unstable manifolds with widely different expansion rates. Unfortunately there are no known examples of geodesic flows that are suitable bases for the perturbation argument; see the discussion at the end of section 4 in \cite{Has-Wil1999}. Along all geodesics in complex hyperbolic space one has Lyapunov exponents of 1 and 2 corresponding to parallel families of planes with curvature $-1$ and $-4$ respectively.
  
  \begin{qst} Is there a metric of negative curvature for which the ratio of largest positive Lyapunov exponent to smallest positive Lyapunov exponent is greater than $2$ for typical geodesics (not just exceptional closed geodesics)?
  \end{qst}
  
  The Anosov splitting  is $C^\infty$ only if and only the metric is locally symmetric. This follows from \cite{Be-Fo-La1992} and \cite{Be-Co-Ga1995}. For surfaces Hurder and Katok \cite{Hur-Kat1990} showed that the splitting must be $C^\infty$ if one of the stable or unstable foliations is $C^{1+o(x\log|x|)}$. In higher dimensions it is expected that the splitting must be $C^\infty$ if it is $C^2$, but this question still seems to be open.

\subsection{ How many   closed geodesics  of length $\le T$ exist? (Communicated by Pollicott).}
Let $(M, g)$ be a closed Riemannian manifold with  negative sectional curvatures. It is known \cite{Margulis2004} that the number $N(T)$  of closed geodesics of length $\le T$ grows
 asymptotically  as ${e^{hT}}/{hT}$.

 \begin{con} {\it $N(T)= (1+ O(e^{-\varepsilon T}))  \displaystyle \int_2^{e^{hT}} \frac{du}{\log u} $. } \end{con} 

The conjecture is true for all
metrics such that the stable and unstable bundles   $E^s $ and $ E^u$ are $C^1$ smooth, which is the case for surfaces and for 1/4 pinched metrics. Question~\ref{C^1quest} above asks whether there are other examples.

\subsection{Infinitely many simple closed geodesics}

\begin{qst}[Miller \cite{Miller2001}, Reid]  {\it Are there infinitely many simple closed 
 geodesics in every hyperbolic complete  three-manifold of finite volume?}\end{qst}

 The answer is positive  for surfaces and estimates for the growth of the number of simple closed geodesics with length $\le L$ are   in  \cite{Rivin2001,Mirzakhani2008}.
   One can ask the same question for manifolds of variable negative curvature in any dimension.
   The answer is  positive for
  generic metrics of negative curvature.

 Reid \cite{Reid1993} has   examples of (arithmetic) hyperbolic manifolds of finite volume in which every closed geodesic is simple.

\section{Non-positive curvature and non-uniform hyperbolicity of the geodesic flow}

Recall that the rank of a vector $v$  in such a manifold is the dimension of the space of Jacobi fields along the geodesic $\gamma_v$ that are covariantly constant. It is easily seen that rank is upper semi continuous. All vectors have rank $\geq 1$, since the velocity vector field is a covariantly constant Jacobi field. The rank of the manifold is the minimum rank of a vector. The set of vectors of minimum rank is obviously open and is known to be dense \cite{Ballmann1982}. The definitions generalize the classical notion of rank for locally symmetric spaces 
of non-compact type.

\subsection{Ergodicity of geodesic flows}

\begin{qst}  \label{ergquest}{\it Is the geodesic flow of a metric of
non-positive curvature on a closed
surface  of genus $\ge 2$ ergodic with respect to the Liouville measure? } \end{qst} 

It is known
that the flow is  ergodic on the  set of geodesics passing through points where the curvature is negative \cite{Pesin1977}.   The complement of this set consists of vectors tangent to zero curvature geodesics, i.e.
geodesics along which the Gaussian curvature is always zero. It is not known whether this  set must have measure
zero. There is an analogous question in higher dimensions: \emph{is the geodesic flow of a closed rank one manifold of non-positive curvature  ergodic}? Here it is known that the flow is ergodic on the set of rank one vectors \cite{Bal-Bri1982, Burns1983}, but it is not known if the complementary set (of higher rank vectors) must have measure zero.

Several papers published in the 1980s (notably \cite{Burns1983}, \cite{Bal-Bri1982}, and \cite{Ba-Br-Eb1985}) stated that this problem had been solved. These claims were based on an incorrect proof that the set of higher rank vectors must have measure zero.

 The measure considered above is the Liouville measure. There is a (unique) measure of maximal entropy for the geodesic flow of a rank one manifold. It was constructed by Knieper, who proved that it is ergodic \cite{Knieper1998}.

\subsection{Zero curvature geodesics and flat strips}

 \ Consider a closed surface of genus $\ge 2$.   

 \begin{qst}\label{flatgeoquest}[Burns]  {\it Is there a  $C^\infty$   (or at least  $C^k, \ k\ge 3)$ metric 
  for which there exists a non-closed geodesic along which the Gaussian curvature is everywhere zero?} \end{qst}

The question is of particular interest when the metric has non-positive curvature, since a negative answer in this case would give a positive answer to Question~\ref{ergquest}.

If we assume that the curvature is only  $C^0$, then a  metric with
 such  a geodesic can be constructed. Each end of the geodesic is asymptotic to a closed geodesic. But a geodesic along which the curvature is zero cannot spiral into a closed geodesic  if the curvature is $C^2$; see \cite{Ruggiero1998}. Wu \cite{Wu2013} has recently given a negative answer to the question under the assumption that the subset of the surface where the curvature is negative has finitely many components; he also  shows under this hypothesis that only finitely many free homotopy classes can contain zero curvature closed geodesics.

In a translation surface, the number of free homotopy classes that contain closed geodesics with length $\leq L$ grows quadratically with $L$. 

 \begin{qst} [Knieper]
 {\it  
Consider a closed surface of genus $\geq 2$ with a metric of non-positive curvature. Let $Z(L)$ be the number of free homotopy classes containing 
closed geodesics with length $\leq L$ along which the curvature is everywhere zero.
Is there a quadratic upper bound for $Z(L)$, i.e.\  is there a constant $C \geq 0$ such that
         $$ Z(L) \leq C \cdot L^2?  $$
   } \end{qst}

A negative answer to Question~\ref{flatgeoquest} in the case of non-positive curvature would imply that  $Z(L)$ is bounded.

A flat strip is a totally geodesic isometric immersion of the Riemannian product of $\Real$ with an interval. Zero curvature geodesics can be viewed as flat strips with infinitesimal width. It is known that all flat strips (with positive width) in compact surfaces with non-positive curvature must close up in the $\Real$-direction; in other words they are really immersions of the product of $S^1$ with an interval. Furthermore for each $\delta > 0$ there are only finitely many flat strips with width $\geq \delta$. Proofs can be found in a preprint of  Cao and Xavier \cite{Cao-Xav}. The main idea is that if two flat strips of width $\delta$ cross each other at a very shallow angle, then their intersection contains a long rectangle with width close to $2\delta$. 

It is possible that there are only finitely many flat strips. This is true when the set where the curvature is negative has finitely many components \cite{Wu2013}, but the only effort to prove it in general \cite{Rodriguez Hertz2003} was unsuccessful. 

{\bf Remark. }  Question \ref{flatgeoquest} still makes sense even if the curvature of the manifold is not restricted to be non-positive. 
It might generalize also to higher dimensions. What can be said about a geodesic along which all sectional curvatures are $0$? Or less stringently the sectional curvatures of planes containing the tangent vector?

\subsection{Flats in rank one manifolds with non-positive curvature}
A flat is a flat strip of infinite width, in other words a totally geodesic isometric immersion of the Euclidean plane.
Eberlein asked whether a compact rank one manifold that contains a flat must contain a closed flat, i.e.~ a totally geodesic isometric immersion of a flat torus. Bangert and Schroeder gave an affirmative answer for real analytic metrics
 \cite{Ban-Sch1991}. 
The problem is open for $C^\infty$ metrics.

\subsection{Besson-Courtois-Gallot}

Does their rigidity theorem from \cite{Be-Co-Ga1995}, Theorem~\ref{BCGthm} in this paper, generalize to  higher rank symmetric spaces of non-compact type?

Connell and Farb \cite{Con-Far2003a, Con-Far2003b} extended the barycenter method, which plays a vital role in  \cite{Be-Co-Ga1995}, and proved that the theorem holds for a product in which each factor is a symmetric space with negative curvature and dimension $\geq 3$. We also refer the reader to their extensive survey in \cite{Con-Far2003c}.

A recent preprint of Merlin \cite{Merlin2016} uses a calibration argument to show that $h_{vol}^4$Vol is minimized by the locally symmetric metric on a compact quotient of the product of two hyperbolic planes.

\section{ Manifolds without conjugate points (rigidity conjectures) }\label{sec8}

There are many rigidity results for manifolds with non-positive curvature that might extend to manifolds with no conjugate points. One still has the basic setting of a universal cover homeomorphic to $\Real^n$ in which any two points are joined by a unique geodesic. However the convexity of the length of Jacobi fields, which is the basis of many arguments used in non-positive curvature, is no longer available. 

Two major results of this nature are:

\begin{thm}\label{Hopf}
A Riemannian metric with no conjugate points on a torus is flat.
\end{thm}

\begin{thm}\label{BangEmm}
 Let $g$ be a complete Riemannian metric without conjugate points on the plane $\Real^2$. Then for every point $p$
 $$
 \liminf_{r \to \infty} \frac{\text{\rm  area}\,B(p,r)}{\pi r^2} \geq 1,
 $$
 with equality if and only if $g$ is flat.
\end{thm}

Theorems~\ref{Hopf} and \ref{BangEmm} are easy for  manifolds with non-positive curvature but require subtle arguments in the context of no conjugate points.

The two dimensional case of Theorem~\ref{Hopf} was proved by E.~Hopf  \cite{Hopf1948} and the general case by Burago and Ivanov 
\cite{Bur-Iva1994}. The Lorentzian analogue of Theorem~\ref{Hopf} is false; two dimensional counterexamples are constructed in 
\cite{Bav-Mou2013}. The Finsler analogue is also false as there are examples of Finsler metrics without conjugate points on the torus that are not the Minkowski metric, see e.g. \cite[\S 33]{Busemann1955}.

 \begin{qst}[Communicated by S. Ivanov] 
 Is the  geodesic flow of  a Finsler metric on the torus   conjugate  to that of the (locally) Minkowski metric on the torus?  
 \end{qst}
 
  All smooth examples we are aware of have geodesic flows   smoothly conjugate to that of the locally Minkowski metric, but generally it is not even known  whether the velocity vectors of    geodesics with irrational rotation numbers  are dense in the tangent bundle.

Theorem~\ref{BangEmm} is a recent result of Bangert and Emmerich \cite{Ban-Emm2013}, which greatly improved on earlier results in \cite{Bur-Kni1991} and \cite{Ban-Emm2011}. 

Bangert and Emmerich's work was motivated by a question, asked by Bangert and
 Burns-Knieper, which is answered by the following theorem, whose proof is based on
Theorem~\ref{BangEmm}.

\begin{thm}[\cite{Ban-Emm2013}]  {\it Consider a complete  Riemannian  metric without conjugate points on the cylinder $\Real \times S^1$. Assume that the ends spread sublinearly,  i.e.\
 $$
 \lim_{d(p, p_0)\to \infty } \frac{l(p)}{dist(p, p_0)} = 0,
 $$
 where $l(p)$ is the length of the shortest  geodesic  loop  based at $p$. Then the metric is flat.}
 \end{thm}

\subsection{ Divergence of geodesics }

 Let $\alpha$ and $\beta$ be two geodesics in a complete simply connected Riemannian manifold without conjugate points.

\begin{qst} {\it Suppose $\alpha(0)=\beta(0)$. Does $dist(\alpha(t), \beta(t)) \to \infty$ as $t\to \infty$?}\end{qst} 

The answer is positive in dimension 2  \cite{Green1954}. The question is open in higher  dimensions; the proof in \cite{Green1956} is incorrect. The answer to the corresponding question about Finsler metrics is negative. The 
Teichm\"uller metric on Teichm\"uller space has no conjugate points, but Masur showed that there exist two rays through the same point which do not diverge \cite{Masur1975}.

A closely related question is:

\begin{qst} {\it Suppose $dist(\alpha(t),\beta(t)) \to 0$ as $t \to -\infty$. Does $dist(\alpha(t), \beta(t)) \to \infty$ as $t\to \infty$?}\end{qst}

Without positive anwers to these questions there would seem to be little hope for a satisfactory analogue of the sphere at infinity, which plays a prominent role in the case of non-positive curvature.

\subsection{Parallel postulate questions}

\

\
The flat strip theorem is another basic tool in the study of manifolds with non-positive curvature.
It states that ``parallel geodesics'' must bound a flat strip. More precisely if
if $\alpha$ and $\beta$ are two geodesics in a simply connected manifold with non-positive curvature such that $dist(\alpha(t),\beta(t))$ is  bounded for all $t \in \Real$, then the two geodesics are the edges of a  totally geodesic isometric immersion of the Riemannian product of $\Real$ with an interval.

The flat strip theorem fails for manifolds with no conjugate points (although it does extend fairly easily  to manifolds without focal points). Counter examples have been given in ~\cite{Burns1992} and by Kleiner (unpublished).  Kleiner's example has a copy of $\mathbb{ Z} \times \mathbb{ Z}$ in its fundamental group, but does not contain a corresponding flat torus. He  perturbs  the Heintze example (Example 4 in \cite{Ba-Br-Eb1985}), which we described in subsection ~\ref{relationsbetween}. 
The following basic question appears to be  open in general (although Proposition 4 in  \cite{Eschenburg1977} suggests an affirmative answer under some extra hypotheses).

\begin{qst}
{\it Must homotopic closed geodesics in a manifold  of dimension $\ge 3$  with no conjugate points be homotopic through closed geodesics?}
\end{qst}
 In dimension 2, the answer is positive. This follows from  a  modification of the Birkhoff minimax procedure recalled  in \S  \ref{mm} 
  combined with the  results/methods of \cite{Morse1924}, as we now explain.   

 Consider  two homotopic  closed geodesics $\gamma_0$ and $\gamma_1$. They have the same length  by   \cite[Theorem 12]{Morse1924}, which we denote by $L_0$
  and their 
 lifts  to the universal cover bound a periodic strip $S$ by \cite[\S 19]{Morse1924}. Next, take 
  a  continuous bijection  $h:S^1 \times [0,1] \to M$ such that $S^1  \times \{0\} = \gamma_0$ and $h(S^1  \times \{1\}) = \gamma_1$.  Then  apply to each curve $h(S^1  \times \{t\})$ a shortening procedure; in dimension 2 one can organize the shortening procedure such that the lift of the 
   projection to the universal cover does not leave the strip $S$. 
   
   The natural modification of the Birkhoff minimax procedure proves then the existence of a closed geodesic homotopic to $\gamma_0$   whose lift to the universal cover lies in    the strip $S$; this geodesic  also  has length $L_0$. 
   Arguing as in \cite[\S 19]{Morse1924}, we obtain that the lift of this geodesic cuts the strip in two  substrips, $S_1$ and $S_2$. Recursively applying  the same arguments, we obtain that the  strip $S$ has a everywhere dense set of nonintersecting geodesics such that any of them cuts the strip into two substrips. Passing to the limit we get a  foliation of the strip by lifts of closed geodesics and hence  the desired homotopy of $\gamma_0$ to $\gamma_1$  through closed geodesics.

Note that for  dimensions $\ge 3$, part of the argument  still works (in particular homotopic closed geodesics in a closed manifold with no conjugate points must have the same length)  and one can show that $\gamma_0$ and $\gamma_1$ belong to the same connected component of the
set of closed curves  of length equal to $L_0$, but it is not clear why the set of closed geodesics homotopic to $\gamma_0$ is path-connected.

One can still hope for  a version of the flat strip theorem in manifolds with no conjugate points in which  there is a large  enough family of ``parallel'' geodesics. One would then hope to find a totally geodesic isometric immersion of the Euclidean plane. Rigidity might hold in the large, even though the examples above show that it breaks down locally.

The simplest  question of this type asks if a Riemannian  plane satisfying Euclid's  5th  postulate must be flat.
Consider the plane $\Real^2$ with a complete Riemannian metric. Assume that for every geodesic and for every point not on the geodesic there exists precisely one non-trivial  geodesic that passes through the point and does not intersect the geodesic.
This assumption   implies  that the metric has no conjugate points.

\begin{qst}[\cite{Bur-Kni1991}] {\it Must a metric satisfying this version of the parallel postulate be  flat?} \end{qst}

The question looks like a question in the synthetic geometry, but it is not, since we do not require a priori that the other axioms of the Euclidean geometry are fulfilled (for example the  congruence axioms).

\subsubsection{Higher rank rigidity}
 \begin{con}[Spatzier]
{\it Let $(M,g)$ be a closed symmetric space of non-compact type and of higher rank.  Then the only metrics on $M$ with no conjugate points are homothetic rescalings of $g$.} \end{con}

\subsection{ Is $M\times S$ with no conjugate points a direct product? }

Consider the product of a closed surface $M^2$ of genus $\ge 2$ with the circle $S$. Let $g$ be a Riemannian metric on $M\times S$ with  no conjugate points.

\begin{con}[Burago-Kleiner] {\it The metric  on the  $\mathbb{Z}$-cover corresponding to $S^1$ factor is  a direct product: $(M\times \Real, g)=(M, g_1)\times  (\Real, dt^2)$.  }\end{con}

One can of course make more general conjectures, e.g.\ about metrics
without conjugate points on products, or on manifolds which admit
non-positively curved metrics with higher rank, but the  above conjecture with 
$M \times  S^1$ seems to be the easiest one (and is probably still hard to prove).

\subsection{ Magnetic geodesics without conjugate points }

One can generalize  the notion ``conjugate points" for magnetic geodesics (even for arbitrary natural  Hamiltonian systems.)

\begin{qst}[Paternain]  {\it Consider a magnetic  flow on  a closed surface of genus $\ge 2$ and an   energy level such that the topological entropy vanishes. Suppose that the magnetic  geodesics    on this level do not have conjugate points. Is the system   locally symmetric in the sense  that the metric has constant curvature and the magnetic form
is a constant multiple of the volume form?} \end{qst}

 Many of the earlier questions in this section   can also be asked about magnetic geodesics lying on a certain (possibly sufficiently high) energy level. 
 Let us mention  \cite{Bialy2000}, 
  where a  natural analog of Theorem \ref{Hopf} was proved under the assumption that the metric is conformally flat, and  it was  conjectured that this assumption is not essential.

\section{Beyond no conjugate points}

\subsection{Ergodic geodesic flows}

{\it Does every closed manifold (with dimension $\geq 2$) admit a metric (Riemannian or Finsler) with ergodic geodesic flow?}

This is known for surfaces~\cite{Donnay1988II}, for $3$-manifolds~\cite{Katok1994}, for product manifolds in which the factors have dimenson $\leq 3$ \cite{Bur-Ger1994}, and for spheres \cite{Bur-Ged}. Donnay and Pugh have shown that  any embedded surface can be perturbed, in the $C^0$ topology to an embedded surface whose geodesic flow is ergodic \cite{Don-Pugh2004}. These  constructions all make essential use of the focusing caps introduced by Donnay in \cite{Donnay1988I}.

The general problem is still wide open despite some reports of its solution (p.~87 of \cite{Berger2000} and Section 10.9 of \cite{Berger2003}).

\begin{qst} {\it \label{q12} 
Is there a Riemannian metric on a closed manifold  with everywhere positive sectional curvatures and ergodic geodesic flow?} \end{qst} 

Metrics close to the standard metric on $S^2$ would be especially interesting. The best result in this direction is due to Burago and Ivanov \cite{Bur-Iva2016}. They exhibit smooth Finsler metrics with positive measure theoretic entropy that are small perturbations of the round metric on $S^4$.

\subsection{ The measure of transitive and recurrent sets  }

Let $(M^2 , g)$ be a closed surface with ergodic geodesic flow. Then every tangent vector lies in
one of the following sets: 

\begin{align*}
T_b &:=\{ v\in SM\mid \textrm{ any lift of $\gamma_v$ stays in a  bounded subset  of $\widetilde M$}\}, \\
T_p &:=\{ v\in SM\mid \textrm{ any lift of $\gamma_v$
 is unbounded, and  approaches  infinity properly}\}, \\
T_i &:=TM\setminus (T_b \cup T_p).
\end{align*}

 All three sets are measurable and invariant under the geodesic flow.  By ergodicity one of them has full measure. 

\begin{qst}[Benjamin Schmidt]  {\it Which of these sets  has  full measure?} \end{qst}

The expected  answer depends on the genus of the surface:
 on the torus the set $T_{i}$ should have  full measure, and on the surfaces of higher genus
 $T_{p}$ should have  full measure.

  In this context it is natural to consider only minimizing geodesics: if one shows that the set of minimizing geodesics has non-zero  measure, then it has  full measure by ergodicity of the geodesic flow.
However, there exist examples of  metrics on the torus such that the set of minimising geodesics has small Hausdorff dimension.

The question seems to make sense in any dimension.

\subsection{Generic metrics} \label{Generic}

We  denote   by $\mathcal G^k(M)$ the space of all smooth 
 Riemannian metrics on $M$ with the  $C^k$-topology.  Riemannian   metrics are sections of the bundle 
$S^2T^*M$. Our topology is     the weak $C^k$-topology on the space of embeddings of $M$ to $S^2T^*M$, see e.g.\ \cite[Chapter II]{Hirsch1976}. The definition can be extended to metrics of all signatures and to Finsler metrics. 

We use $C^k$-topology either locally or on closed manifolds, where we may 
assume that the atlas is finite and the gluing functions have all derivatives bounded.   Then the   convergence in the  topology means the following: a sequence $g^{(1)},g^{(2)},\dots$ of metrics converges to a metric $g$ in the topology $C^k$ if  in  local charts  the entries $g^{(m)}_{ij}$ (viewed as local functions in this chart) and their partial derivatives  $\frac{\partial g^{(m)}_{ij}}{\partial x_\ell}, 
\frac{\partial^2 g^{(m)}_{ij}}{\partial x_\ell \partial x_s},\ldots $ up to order $k$ converge to those of $g_{ij}$.

  \begin{qst} Assume $M$ is closed. 
  Do the  metrics with positive entropy  form a  $C^k$-dense subset of $\mathcal G(M)$? 
  \end{qst}

  This question can be asked for any $k$. For $k=2$ it was recently positively answered in all dimensions in 
  \cite{Contreras2010}; see also \cite{Con-Pat2002}. For $k = \infty$ and dimension 2, it was answered positively in  \cite{Kni-Wei2002}.

\begin{qst}  
Assume $M$ is closed. 
If $\dim(M)\ge 3$,  do  the  metrics for which  the  geodesic flow is  
transitive  (i.e.\ has an orbit whose closure is $SM$)   form a  dense subset of  $\mathcal G(M)$?
\end{qst} 

 This question is  closely related to the famous 
 Arnold diffusion conjecture.
If $\dim(M)=2$ and $k$ is large enough, 
 the answer is negative: if we take a metric with  integrable geodesic flow  such that certain 
  Liouville tori are irrational,  and 
 any small perturbation of the metric has non-transitive geodesic flow   by the KAM theory.

\begin{qst} {\it  Assume $M$ is closed.
  Is  there a dense  subset of metrics in $\mathcal G$  for which  the tangent vectors to  the closed geodesics are dense in $SM$? }
  \end{qst} 

 It is known that  this is the case for metrics with non-uniformly hyperbolic geodesic flows.  
Note that is relatively easy to construct, see e.g.\ \cite{Weinstein1970},  a metric on any manifold such that a certain non-empty open subset of $SM$ contains no vectors tangent to a closed geodesic.

\subsection{Density in the manifold}

In the previous section we discussed the unit tangent bundle $SM$. In this section we ask similar questions  about $M$ itself, and we will not assume that the metrics are generic.

\begin{qst} {\it  
Is the union of the closed geodesics always dense in the manifold?}\end{qst}

\begin{qst} {\it 
Does every metric on a compact surface of positive genus have a geodesic that is dense in the surface?}
\end{qst} 

\subsection{Finsler metrics with positive flag curvature}

Flag curvature is the Finslerian analogue of sectional curvature. It was shown in \cite{Akbar-Zadeh1988} that: any closed Finsler  surface with constant flag curvature zero is locally Minkowskian; any sufficiently smooth closed Finsler surface with constant negative curvature is Riemannian (there are non-Riemannian examples which are not sufficiently  smooth);  and any  connected closed Finsler surface with constant positive flag curvature   is homeomorphic to $S^2$ or $\Real P^2$. 
The result for negative flag curvature also follows from \cite{Foulon1997} and is true in all dimensions.

\begin{problem} {\it 
Describe  Finsler metrics on $S^2$ of  flag curvature $1$.}
\end{problem}
 
 If such  a metric is reversible, then it is the round Riemannian metric \cite{Bryant2006}.
The Katok examples described in Section~\ref{exkatok} have constant flag curvature $1$.
All projectively flat metrics of constant flag curvature on $S^2$ are described in \cite{Bryant1997}. Bryant has also constructed examples that are not  projectively flat
\cite{Bryant1996,Bryant2002}. Unlike Katok's example, all of the examples in Bryant's papers are Zoll, i.e.\ all of their closed geodesics are closed and of the same length. 

Even though the Finsler metrics of constant flag curvature $1$ on $S^2$ are not fully classified, the possible behavior of their geodesic flows is understood. It is shown 
in \cite{BFIMZ2019} that the geodesic flow of a metric of constant curvature $1$ is  conjugated to that of the Katok metric. This implies that 
 one of the following two possibilities must occur:
  \begin{itemize}  \item[(i)]  Precisely two geodesics are closed. These two  geodesics are a simple closed curve parametrized in opposite directions. This curve is tangent to a  non-trivial Killing vector field which vanishes precisely at two points.  
 All of the non-closed geodesics have the same asymptotic rotation number $\rho$ around the fixed points of the Killing field. There are distinct irrational numbers $\alpha_1$ and $\alpha_2$  such that the lengths of the two closed geodesics are  $\alpha_1\pi$ and $\alpha_2\pi$.

    \item[(ii)]  All geodesics are closed. There are  an even  integer
     $k$  and positive rational numbers $\beta_1$ and $\beta_2$  with $\tfrac{1}{\beta_1} + \tfrac{1}{\beta_2}=1$
     such that 
at most   two   geodesics have  length different from  $k\pi$ and the lengths of the exceptional geodesics are $\beta_1\pi$ and $\beta_2\pi$.   
In the case   $\beta_1=\beta_2=2$, all geodesics are closed and have the same length $2\pi$; we will call such metrics {\it Zoll metrics}.  
\end{itemize} 

Moreover, the Finsler metrics satisfying (i) necessary admit a Killing vector field and are quite well understood. Applying the Zermelo transformation  to such metrics (see e.g. \cite{Fou-Mat2018}), one  obtains non-Zoll metrics satisfying (ii). Moreover, any such  metric  (of flag curvature $1$ on $S^2$  and satisfying (ii)) admitting a Killing vector field  can be obtained    by this procedure.  Unfortunately, no other  examples  of non-Zoll metrics satisfying (ii) are known.

\begin{problem}
{\it Construct   examples  of Finsler metrics of constant curvature on $S^2$ that are not Zoll and such that they admit no Killing vector field. }
\end{problem}

For recent progress in this direction see \cite{Lan-Met2018}.

In a closed (even complete) Riemannian manifold  $M$ of positive sectional curvature two 
closed totally geodesic submanifolds  such that the sum of their dimensions is at least the dimension of $M$ must intersect  \cite{Frankel1961}.

\begin{qst} {\it \label{Fra}  Does Frankel's Theorem remain  true for Finsler manifolds? More precisely,  must two closed  
 totally geodesic submanifolds of a  closed Finsler manifold
 with positive flag  curvature   intersect if the sum of their dimensions is at least the dimension of the ambient manifold? } 
\end{qst}

In dimension $2$ the result was established in \cite{BFIMZ2019}.  
A positive answer to this question was claimed in  \cite{Koz-Pet2000}, but   the proof  has a serious flaw, which is explained in \cite{BFIMZ2019}. Note that there are 
Finsler surfaces with positive flag curvature that contain non-intersecting closed 
geodesics, see \cite{BFIMZ2019, Rademacher2016} This is possible because a curve must have geodesic parameterizations in \emph{both} directions in order to be totally geodesic,
  non-reversible Finsler metrics can therefore have simple 
closed geodesics that are not totally geodesic submanifolds.

The simplest still open  version of Question \ref{Fra} is the following question. 

\begin{qst} {\it  Let $(M^3, F)$ be a closed 3-dimensional reversible  Finsler manifold of positive flag curvature and $N$ be a closed totally geodesic 2-dimensional submanifold. Is it true that each closed geodesic intersects this submanifold or lies on it?
 } 
\end{qst}

\subsection{ Gaidukov in higher dimensions.}

We consider a
 Riemannian metric on an oriented  closed  surface $M$  of genus $\ge 1$. Let $\Gamma$ be a non-trivial
 free homotopy     class and $p$ a point in $M$. A  theorem of Gaidukov  \cite{Gaidukov1966} says that
 there exist a closed geodesic $\gamma:\Real\to M$ in $\Gamma$ and a ray $\beta:[0,\infty)\to M$
with $\beta(0)=p$  such that $dist(\beta(t), \gamma(t)) \to 0$ as $t \to \infty$.   As explained in \cite{Bia-Pol1986}, Gajdukov's results follow 
 easily from the classical results of \cite{Morse1924} and  \cite{Hedlund1932}.

\begin{problem}[Benjamin Schmidt] \label{Sch} {\it Generalize this statement to higher dimensions.} \end{problem}

 Gaidukov's proof is profoundly two-dimensional and cannot be generalized.  But investigations in the direction of Problem \ref{Sch} are important and
plentiful in Aubry-Mather theory and Arnold diffusion.     Instead of a single closed geodesic one should consider a minimizing set, namely  the support of  a  minimal  measure.
See e.g.\ \cite{Mather1991,Con-Itu1999,Mather2004, Bernard2010}. 
 
 \subsection{ Systolic and diastolic inequalities  for  surfaces  (communicated by Guth,  Rotman and Sabourau)}  
 
 Let $(M,g)$ be a compact Riemannian surface. The inequalities in question compare the length of certain short closed geodesics with the square root of the area of $(M,g)$,
 Recall that the  systole $sys(M,g)$   is 
  the least length of a non-trivial closed geodesic. Two other geometrically meaningful constants can be defined by the following minimax procedure:
  $$
L(M,g) = \inf_f \max_{t\in\Real} F[f^{-1}(t)],
$$
in which the infimum is taken over all (Morse) functions $f: M \to \Real$ and the functional $F$ is either (a) the total length of $f^{-1}(t)$ or (b) the length of its longest component. In both cases $L(M,g)$ is realized as the length of a certain union of closed geodesics. In case (a), $L(M,g)$ is one definition of the diastole $\textrm{dias}(M,g)$ of the Riemannian surface (at least two different notions of diastole have been studied; see \cite{Bal-Sab2010}).

By a result originally due to   \cite{Croke1988} and later improved in  \cite{Nab-Rot2002}, \cite{Sabourau2004} and \cite{Rotman2006}, for
every Riemannian metric $g$ on the sphere $S^2$ one has
$$
 \textrm{sys}(S^2, g) \leq \sqrt{ 32} \sqrt{\textrm{area}(S^2, g)}. 
$$  Actually, in \cite{Croke1988} it was suggested that the  constant  $\sqrt{32}$ in the above  inequality  can be replaced by   
 $\sqrt{2\sqrt{3}}$. The following example due to \cite{Croke1988}  shows that one can not go below $\sqrt{2\sqrt{3}}$: take two congruent equilateral triangles and glue them along their  boundaries. 
The example   is not smooth and suggests the study of  systolic inequalities on the space of smooth metrics with conical singularities. On the space of such metrics, let us  consider the function 
$$
\sigma(g)= \frac{\textrm{area}(S^2, g)}{ \textrm{sys}(S^2, g)^2}. 
$$
This function has many nice properties; for example it is Lipschitz with respect to the appropriate distance on the space of metrics.  By  \cite{Croke1988}, the function has a positive minimum, and the natural conjecture is that the minimum 
is attained on the sphere constructed from two equilateral triangles as  described above.

Other critical points of this function are also interesting: Balacheff showed that the round metric on the sphere
 is a critical point  \cite{Balacheff2006}.  Moreover, by \cite{Alv-Bal2014}, Zoll Finsler manifolds $(M,  F)$  of all dimension  satisfy  the following alternative: for any  Finsler deformation $F_t$  with $F_0=F$   either the function  $ t \mapsto \sigma(F_t)$   has a
strict local minimum  at $t = 0$, or  the family $F_t$   is tangent up to every order to the space of Zoll metrics. In dimension 2, the result is further improved in \cite{ABHS2017,ABHS2018}, where it is shown that both  in the Riemannian  and in the Finslerian cases  the function $\sigma(F_t)$ defined above    has a strict local minimum  at $t = 0$ if $F$ is a Zoll metric.

Note though that in the papers \cite{Alv-Bal2014,ABHS2017,ABHS2018} the metrics  and the deformation are assumed to be smooth (at least $C^3$). In view of the Croke's example, it may be natural to consider also non-smooth metrics and non-smooth deformations.    

\begin{qst}[Babenko] {\it
 Does there exist,  in the space of metrics with conical singularities, a (continuous or smooth) family of 
  metrics $g_t$ such that $g_0$ is the round metric of the sphere, $g_1$ is the metric from the example above,  and 
   $\sigma(g_t)$ is a decreasing function of $t$.}
\end{qst}

A simpler version of this question is

\begin{qst}[Babenko] { \it  Does every Riemannian metric on $S^2$ with area $4\pi$ that is close enough to the round metric have  a closed geodesic with length $\leq 2\pi$?  In other words, is it true that $\sigma \ge  1/\pi$ for metrics close to  the round metric?  \label{babenko1}}
\end{qst}

The answer depends on what is meant by ``close''.
If ``close'' means $C^2$-close, then the positive answer follows from  \cite{ABHS2017}. In that paper  the condition ``close''  is 
given explicitly in terms of a pinching condition on the curvature. It is proved that  in the neighborhood of a Zoll metric consisting  of  metrics such that 
this  pinching condition is satisfied,  the equality $\sigma =  1/\pi$
 holds if and only if the metric is Zoll.

A corollary of one of the main results in \cite{ABHS2018} generalizes this to Finsler metrics and a $C^3$-neighborhood of the Zoll metrics : more precisely, $\sigma \ge  1/\pi$ in a small $C^3$-neighbourhood of all Zoll  Finsler metrics on $S^2$, with equality if and only if the metric is Zoll.

So a related question is as follows: under which topologies is the round metric  a local minimizer of $\sigma(F_t)$ or $\sigma(g_t)$? 
Does this remain true for weaker topologies than $C^2$, e.g.\ $C^1$ or Gromov-Hausdorff?

For surfaces of higher genus, there exists a constant $C$ such that
\begin{equation}  \label{gr}
\textrm{sys}(M,g) \leq C \frac{\log( \textrm{genus}(M))}{ \sqrt{\textrm{genus}(M)}}{\textrm{area}(M, g)};
\end{equation}
the dependence on the genus of $M$ in this inequality is sharp.  In \cite[Note 7.2.12]{Berger2003} inequality \eqref{gr} is attributed to Gromov and it is explained from what papers, results and ideas of Gromov \eqref{gr} follows.  Section  7.2.1 of \cite{Berger2003} is a very nice survey on systolic inequalities for surfaces.

For the diastole defined by (a), Balacheff and Sabourau showed in \cite{Bal-Sab2010} that there is a constant $C$ such that
$$
\textrm{dias}(M,g) \leq C\sqrt{\textrm{genus}(M)}\sqrt{\textrm{area}(M, g)};
$$
 the dependence of this inequality on the genus is again optimal.
\begin{qst}[Guth] {\it Is there a constant  $C$  such that the invariant defined by (b) above is bounded from above by} $C\sqrt{\textrm{area}(M,g)}${\it ?}
\end{qst}

An affirmative answer to Guth's question would mean that the three quantities under consideration all depend on the genus in different ways, and are therefore measuring different features of the surface. A positive answer would also show that one can
always find a pants decomposition of a closed Riemannian surface of genus
$g$ into $3g-3$ disjoint closed geodesics of length at most $C \sqrt{\textrm{area}(M)}$.
This would precisely give the optimal Bers' constant for a genus $g$
surface. Even for hyperbolic surfaces, this question is still open (see \cite{Buser1992} 
 for partial results).

 The question  above also makes sense in higher dimensions. In this case $f^{-1}(t)$ would be 
an $(n-1)$-complex, where $n$ is the
dimension of the manifold, and $F$ would measure its total volume or the volume of its largest component.

\subsection{Questions related to the systole   in higher dimensions}

\begin{qst}[Question 4.11 of \cite{Gromov2001}] { \it Does a Riemannian metric on a real projective space with the same volume as the canonical metric have  a closed geodesic with length $\leq \pi$?}
\end{qst}

In dimension two, an affirmative answer is in \cite[Proposition 4.10]{Gromov2001}. 

\begin{qst}[\'Alvarez Paiva]  { \it Can there be  Riemannian metrics on $S^3$ or $S^2 \times S^1$ with unit volume all of whose closed geodesics are long? More specifically, does every  metric on these spaces have  a closed geodesic with length $\leq 10^{24}$ if the volume is $1$?}
\end{qst}

\subsection{Metrics such that all geodesics are closed} 

This is a classical topic --- the first examples go back  at least 
 to \cite{Zoll1903}; see \cite{Besse1978} for details. The book   \cite{Besse1978} is still up to date, and many problems/questions  listed in it (in particular in Chapter 0 \S D) are still open. As the 
  most interesting question from their list  we suggest:
  
  \begin{qst} {\it Let $M$ be a closed manifold not covered by  a sphere. Let $g$ be a metric on $M$ for which all geodesics are closed. Is $(M,g)$ a  CROSS (=compact
rank one symmetric) manifold (with the standard metric)?}
\end{qst}

Of course,  variants of this question can be asked about Finsler and Lorentzian metrics. In the Finsler setting, one may ask for example 
 to describe all Finsler metrics on $CP(n)$ such that their geodesics are geodesics of the standard (Fubini-Study)  metric on $CP(n)$.  Of course, the answer in the Finslerian case is expected to be much more difficult; for example by \cite{Bryant1996} on the 2-sphere  there exists  
  a  two-parameter family of projectively flat Finsler metrics  on the sphere, such that the flag curvature is $+1$.

  In the Lorenzian setting, one can ask to construct all manifolds such that all light-like geodesics are closed, see e.g.\ \cite{Mou-Suh2013, Suhr2013b}.

An easier version of the question above is the following conjecture sometimes attributed to Berger: 

\begin{con} {\it Let $(M,g)$ be a closed  simply connected manifold such that  all geodesics are closed. Then all   geodesics have the same length.} 
\end{con}

The conjecture is proved when $M$ is the sphere of dimension $\ne 3$ in \cite{Gro-Gro1981, Rad-Wil2015}.

\section{Lorentzian metrics  and metrics of arbitrary signature.} 

\subsection{ Closed geodesics.}
Most of the questions we asked about the Riemannian and Finsler metric can be modified such that they are also interesting in the 
 semi-Riemannian metrics (i.e.~of arbitrary signature) and in particular when the signature is Lorentzian.  
It appears though that the answers in the Lorentzian case are sometimes very different from those in the Riemannian case. Many methods that were effectively used in the Riemannian case, for example the variational methods, do not work in the case of general  signature. 

 Let us consider as an example the question analogous to the one we considered in Section 2: 
how many geometrically different closed geodesics must there be for  a Lorenzian metric on a closed manifold. 

First let us note that there are two possible natural notions of closed geodesic in the Lorentzian setting: one may define 
a closed  geodesic as an  inmmersion $\gamma:S^1 \to M$ such that  $\nabla_{\dot \gamma} \dot \gamma=0$, or as a curve $\gamma:S^1 \to M$ 
 that can be locally reparameterized in such a way that it satisfies the equation $\nabla_{\dot \gamma} \dot \gamma=0$. 
 
 In the Riemannian case, these definitions are essentially equivalent, since this reparameterization can always be made ``global''. If the signature is indefinite, it is easy to construct 
  examples of 
 an embedding  $\gamma:S^1 \to M$   such that locally $\gamma$ can be reparameterized so that it becomes  an affinely-parameterized  geodesic, but globally such reparameterization is impossible: 
  the velocity vector of the geodesic after returning to the same point is proportional but not equal to the initial velocity vector. Of course, this is possible only if the velocity vector is light-like.

We will follow most publications and define a closed geodesic as   an  embedding $\gamma:S^1 \to M$ such that  $\nabla_{\dot \gamma} \dot \gamma=0$. 

To the best of our knowledge, the existence of a closed geodesic on a Lorentzian manifold is a quite complicated problem and nothing is known in dimensions $\ge 3$. In dimension 2, a closed orientable manifold with a metric of signature $(+,-)$ 
 is homeomorphic to the torus. By \cite{Suhr2013a}, every Lorentzian 2-dimensional torus has at least two simple  closed geodesics one of
which is definite, i.e. timelike or spacelike. Moreover, explicit examples show   the optimality of this claim.

We therefore ask the following  
\begin{qst}  
{\it Does every closed  Lorentzian manifold  have at least one closed geodesic?}\end{qst} 

One of course can ask this question about manifolds of arbitrary signature and also about complete (note that in the Lorentzian setting there are many different non-equivalent notions of completeness) manifolds of finite volume. 

See the introduction to the paper \cite{Fl-Ja-Pi2011} for a list of known results about the existence of closed time-like geodesics under additional assumptions.

\subsection{ Light-like geodesics of semi-Riemannian metrics }

Let $M$ be a closed manifold with   semi-Riemannian  metric $g$ of indefinite signature.

\begin{qst} \label{qst:ab}  {\it Can there exist a  complete  semi-Riemannian metric $g$ and a non-trivial 1-form $\eta$  on a closed manifold  such 
  that for  every light-like geodesic $\gamma(t)$ 
   the function $\eta(\dot\gamma(t))$ grows  linearly in both directions: i.e.~ for   every light-like geodesic there exist  $C_1 \neq 0, C_2 $ such that
$\eta(\dot\gamma(t))= C_1\cdot t + C_2$?} \end{qst}

A negative answer to  this question would  give an easy proof of the semi-Riemannian version of the projective Lichnerowicz-Obata conjecture:

{\it  Let a connected  Lie group $G$  act on a closed
  manifold  (of dimension at least $2$) by projective
transformations (diffeomorphisms  that take geodesics   to geodesics without necessarily preserving the parameterization).    Then $G$  acts by isometries, or $g$
has  constant  sectional  curvature. }

If $g$ is Riemannian, the conjecture was proved in \cite{Matveev2005,Matveev2007}.  For metrics of Lorentzian signature, the conjecture was proved in dimension 2 in \cite{Matveev2012b}  and for all other dimensions in \cite{Bo-Ma-Ro2015}. 
 The proof is complicated, and does not generalize for other signatures. 
 In the semi-Riemannian case the following argument gives a proof  for closed manifolds provided that  the answer to the Question  \ref{qst:ab} above is positive.  

It is known  (see, for example,  \cite{Matveev2007}) that a 1-form   $\eta_i$
 generates a one-parameter group of projective transformations, if and only if 
 $$\eta_{i,jk}+ \eta_{j,ik} - \frac{2}{(n+1)} \eta^\ell_{\ ,\ell k} g_{ij}   = \lambda_i g_{jk} + \lambda_j g_{ik}$$
 (for a certain 1-form  $\lambda_i$). 
  We take a light line geodesic $\gamma(t) $,  multiply the above equation by  
 $\dot\gamma^i\dot\gamma^j\dot\gamma^k $   (and sum with respect to repeating indexes) at every point $\gamma(t)$ of the geodesic. 
  The terms with the metric $g$ disappear since $g_{ij} \dot\gamma^i\dot\gamma^j=0$, so we obtain the equation 
 $\frac{d^2}{dt^2}\eta(\dot\gamma(t)) =0$ implying  that $\eta(\dot\gamma(t)) = t\cdot C_1  + C_2$.   A negative answer to the question above implies  that $C_1=0$. Hence  $\eta(\dot\gamma(t))$ is constant, which in turn implies that $\eta_{i,j}+ \eta_{j,i}$ is proportional to $g_{ij}$. Then the covector field $\eta_i$ generates  a one-parameter group of conformal
  transformations.  Finally, the proof of the conjecture follows from a classical observation of H. Weyl  \cite{Weyl1921} that 
   every transformation that is projective and conformal is a homothety.

A version  of the question above is  {\it whether, 
 for a complete semi-Riemannian metric on 
a closed manifold $M$,    the tangent vector of almost every geodesic remains in a bounded set of $TM$.  }  A positive answer to this question immediately implies that the answer to the initial question is negative, thereby proving the projective Lichnerowicz-Obata conjecture on closed manifolds.

\subsection{Completeness of closed manifolds of arbitrary signature (communicated by H. Baum).} 
It is well-known that  a closed Riemannian manifold is geodesically complete, in the sense that every geodesic  $\gamma:(a,b)\to {M}$ can be extended to a geodesic  $\tilde \gamma:\mathbb{R} \to  M$ such that $\tilde \gamma|_{(a,b)}= \gamma$. It is also well known that for  any indefinite signature there exist metrics on closed manifolds that are not geodesically complete. 

\begin{qst} {\it  What geometric assumptions imply that a metric (possibly, of a fixed signature) on a closed manifold is geodesically complete?} 
\end{qst}   

We of course are interested in   geometric assumptions that are easy to check or which are fulfilled for many interesting metrics. 

 Let us mention a few classical results.
 For compact homogeneous manifolds, geodesic completeness was established in 
\cite{Marsden1972}. 
For  compact Lorentz  manifolds  of constant curvature 
geodesic completeness was proved in \cite{Carriere1989} (flat case) and \cite{Klingler1996} (general case).

We refer to the paper \cite{Sanchez2015} for a list of  interesting results on this topic and for a list of open questions from which we repeat here only: 

\begin{qst}[\cite{Sanchez2015}] {\it   Assume that a compact Lorentzian manifold is globally conformal to a
manifold of constant curvature. Must it be geodesically complete?} \end{qst} 

Note also that in the non-compact case homogeneous manifolds of indefinite signature are not necessary geodesically complete; see for example 
\cite[Example 2 in \S 4]{Sanchez2015}. It is interesting  to understand whether completeness of a homogeneous manifold  can follow from algebraic properties of the isometry group.

\section{Integrability and ergodicity of geodesic flows on surfaces of higher genus}

\subsection{Metrics  with integrable geodesic flow on surfaces of genus $\ge 2$}

\begin{qst}[Bangert] {\it Does there exist a Riemannian metric on a closed surface of genus $\ge 2$ whose geodesic flow is completely integrable? }\end{qst} 

The answer may depend on the smoothness of the metric and on what we understand by  ``completely integrable": whether the integral is functionally independent of the Hamiltonian on an open everywhere dense subset, or we additionally assume that the subset has  full measure.

Ma\~n\'e showed that
 a Hamiltonian flow on surfaces is generically Anosov or has  zero Liouville exponents \cite{Mane1996}. This   suggests  that  an   easier version of the above   question would be  the following

\begin{qst}[Paternain] {\it Is there a Riemannian metric on  a closed surface of genus $\ge 2$ whose geodesic flow has zero Liouville entropy?} \end{qst}

Bangert and Paternain have outlined some ideas about what a Finsler metric with this property might look like; finding a Riemannian metric is certainly harder.

A special  case of the Finsler version of the  question above would be:

\begin{qst} \label{question:lang}  {\it  Let $F_1$, $F_2$ be Finsler metrics   on a closed surface of genus $\ge 2$. Assume every (unparameterized) $F_1$-geodesic is an $F_2$-geodesic. Must $F_1$ be obtained from $F_2$ by  multiplication by a constant and adding a closed form? } \end{qst} 

One  can also  ask this    
 question  for arbitrary closed manifolds that can carry a hyperbolic metric (if $F_1$, $F_2$  are  Riemannian, the answer is affirmative  \cite{Matveev2003a}).

The previous question is related to the other questions in this section because of the following observation from \cite{Mat-Top1998}: one can use the second metric to construct an  integral  of the geodesic flow of the first one.   If both metrics are  Riemannian,  the integral is quadratic in velocities and the affirmative answer follows  from \cite{Kolokoltsov1983}; see \cite{Mat-Top2000}.

\begin{qst} \label{question:landsberg} {\it Does there exist a non-Riemannian Finsler metric satisfying the Landsberg condition 
on a surface of genus $\ge 2$?} 
\end{qst} 

The Landsberg condition is defined in \cite{Bao2007}. It implies the existence of an integral for the geodesic flow of the metric, which is the relation of this question with the present section.  See \cite{Gom-Rug2013} for a negative
 answer to the question assuming non-existence of conjugate points.   Note    that  the  natural analog of  Question  \ref{question:landsberg}  is completely open in all settings: locally, globally, in all dimensions, and is one of the most important problems in the Finsler geometry, see e.g.\  \cite{Bao2007}.

In the class of  real-analytic metrics,   the papers \cite{Lang2019, Paternain1997}  answer positively  the Questions \ref{question:lang}  and \ref{question:landsberg}.

\subsection{Integrable geodesic flows with integrals of higher degree }

\begin{con}[\cite{Bo-Ko-Fo1995}]   {\it If the geodesic flow of 
a Riemannian metric on the torus $T^2$  admits an
integral   that is  polynomial of degree 3 in the velocities,  then the metric admits  a  Killing vector field.  } \end{con}

A Killing vector field $V$ allows us to construct an  integral 
$$
I:TM \to \Real, \ I(\xi):=g(V, \xi)
$$
  that is evidently linear  in velocities; its third power is then an integral cubic in velocities.

The motivation to study metrics admitting integrals polynomial in velocities comes from the following observation (which  dates back at  least to Darboux and Whittaker): if the geodesic flow admits an integral that is analytic  in velocities, then each  component of this integral that is homogeneous in velocities  is also an integral. The natural idea is then to    study the integrals  that are polynomial in velocities of low degree.

By the result of Kolokoltsov \cite{Kolokoltsov1983}, no metric on a surface of genus $\ge 2$ admits an integral that is polynomial in velocities and  is functionally independent of the energy integral. 
 The state of the art  if the surface is  the sphere or  the torus can be explained by the following table: 

\begin{center}
\begin{tabular} {|c|c|c|}  
\hline  
         & Sphere $S^2$  &Torus $T^2$ \\  
         \hline  
Degree 1 & All is known & All is known \\ \hline  
 Degree 2 & All is  
known & All is known\\ \hline  
 Degree 3 & Series of examples & Partial negative results \\  
\hline Degree 4 & Series of examples & Partial negative results \\  
\hline Degree $\ge$ 5 & Nothing is known  & Nothing is known \\  
\hline  
\end{tabular}

\end{center} 
\vspace{3ex}  
  
In the table, 
``Degree''  means the smallest degree of a non-trivial  integral polynomial in velocities.  
``All is known'' means that there exists an effective description and classification 
(which can be found in \cite{Bo-Ma-Fo1998}).

A simpler version of the question   is when we replace the geodesic flow in the question above by a   Lagrangian  system with  the Lagrangian  of the form $L(x, \xi):=K +  U=\sum g_{ij}\xi^i \xi^j + U(x) $. In this case we assume  that the integral is a  sum of   polynomials  of degrees 3 and 1 in velocities. The
``partial negative results''  in the table above correspond to this case; moreover,   in most cases it is assumed that the metric $g_{ij}$ is flat, see 
 for example   \cite{Bialy1987,Mironov2010,De-Ko-Tr2012} (and \cite{Bialy2010} for results  that do not require this assumption).

Similar  conjectures   could be posed for integrals of every  degree. If the degree $d$
is odd, the conjecture is that the existence of an integral that is polynomial of the degree $d$ in velocities implies the existence of a Killing vector field.  If the degree $d$
is even, the conjecture is that the existence of an integral that is polynomial of the degree $d$ in velocities implies the existence of an integral quadratic in velocities and not proportional to the energy integral.

Note also that by \cite{Kru-Mat2016} a $C^2$-generic metric does not admit a non-trivial integral which is  polynomial    in velocities.

 Further open problems on integrable systems and integrable geodesic flows  can be found in the recent survey \cite{BMMT2018}.

\subsection{ Integrable billiards and the Birkhoff conjecture} 

Let  $\Omega\subset  \mathbb{R}^2$
 be a strictly convex domain bounded by a smooth curve. 
We consider the (internal) billiard problem.
 The ball (a  point) moves with unit speed and no
friction following a rectilinear path inside the domain.  When the ball hits the
boundary, it is reflected elastically according to the law of optical reflection: the
angle of reflection equals the angle of incidence.

We ask  whether there exists a  non-constant smooth 
 function on the tangent bundle which is  constant along every trajectory. It is well known, see e.g.\ \cite{Tabachnikov2005}, that a billiard in an ellipse is integrable. 
 
 \begin{con}  {\it 
If the billiard in $\Omega$ is  integrable, then the boundary of $\Omega$ is an ellipse. }  
 \end{con} 

 The special case  of the conjecture in which  the caustics of the billiard trajectories are assumed to foliate the whole of $\Omega$ except for one point was proved in \cite{Bialy1993}. Recently, 
essential progress is due to \cite{Av-Si-Ka2016}, where the conjecture was proved  for tables
bounded by small perturbations of ellipses of small eccentricity. 

One can also assume that the integral is polynomial in velocities; non-trivial results in this direction include  \cite{Bolotin1990, Bia-Mir2017}.

See also \cite{Gutkin2012} for  generalizations and other problems related to billiards.

\subsection{Metrics such that one can explicitly find all geodesics} 

Geodesics of a metric are solutions of a non-linear  ordinary differential 
equation  $\nabla_{\dot \gamma} \dot \gamma=0$  which can not be {\it explicitly} solved in most cases. There are only a few 
 examples of 2-dimensional metrics for which one can explicitly find {\it all} geodesics using elementary functions: they are metrics of constant curvature and Darboux-superintegrable metrics (see \cite{Br-Ma-Ma2008} for definition). 
 An interesting problem is to construct other examples of metrics such that all geodesics are explicitly known. A related problem  (asked recently in \cite{Tao2010}) 
 is to construct metrics with an explicitly given distance function.

\section{ Stationary  nets (Communicated by Rotman).}

A graph $G$ in a Riemannian surface  $(S,g)$
 is called a \emph{ stationary net,} if every edge is a geodesic and if at every vertex the sum of diverging unit vectors is zero. Vertices must have valence at least $3$.
 
 {\includegraphics[scale=.35]{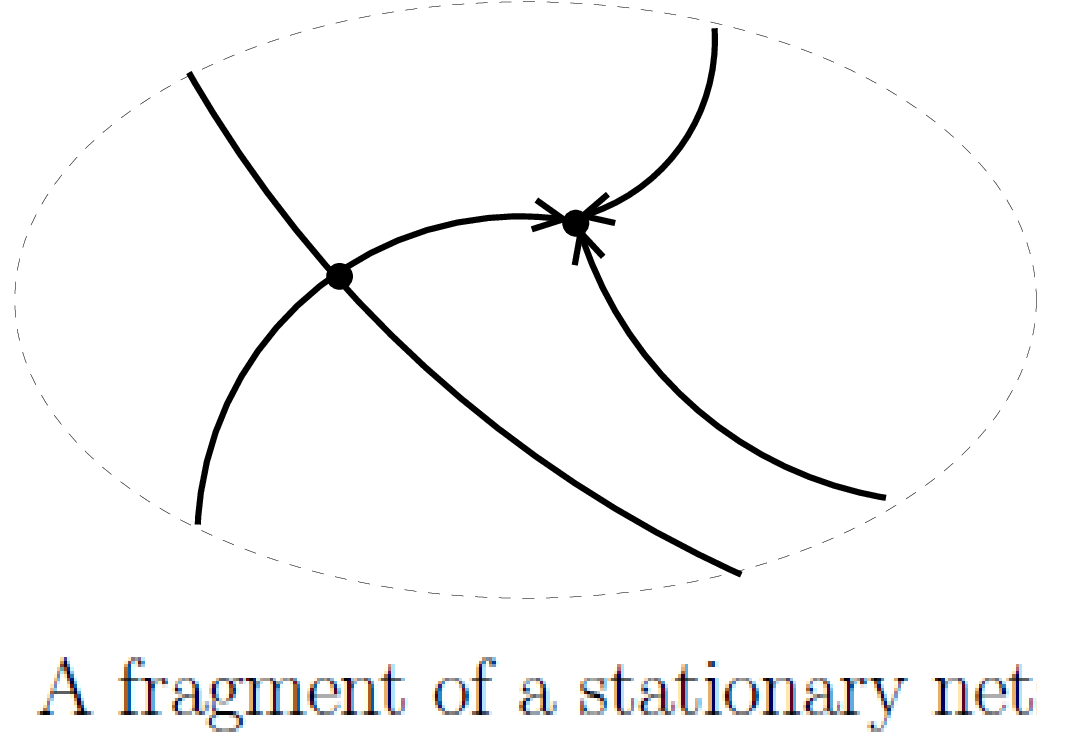}}

\subsection{Stationary $\Theta$-nets. } A $\Theta$-graph is a graph that looks like the Greek
 letter ``Theta'': two vertices connected by three edges, see the picture below.

 \begin{qst}[\cite{Has-Mor1996}]  {\it  Does every metric on $S^2$ admit a stationary $\Theta$-graph?} 
 \end{qst}

Partial results in this direction were obtained in \cite{Has-Mor1996}: they showed  that on every sphere of positive curvature there exists a stationary $\Theta$-net, or a stationary {eight curve} net, or a stationary   eyeglasses net. 

{\includegraphics[scale=.35]{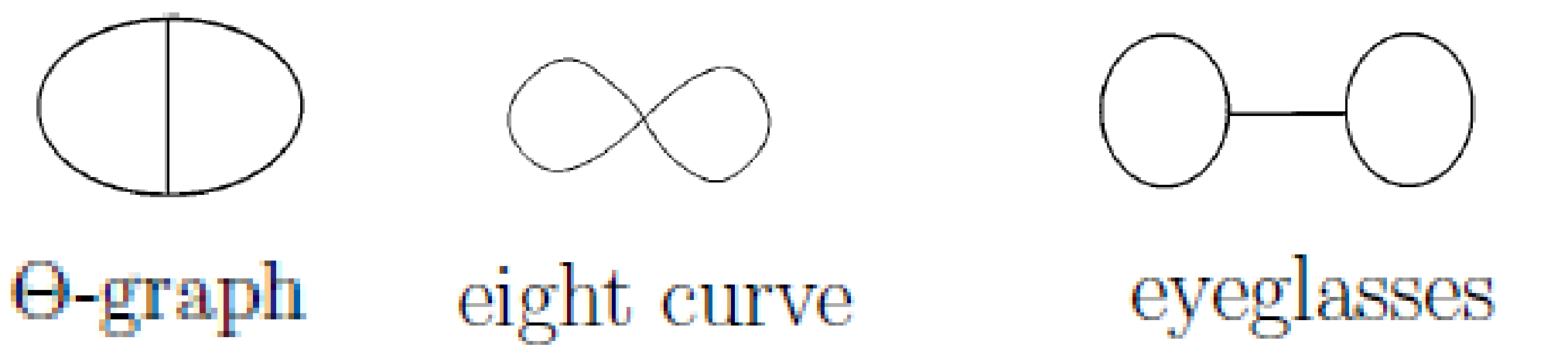}}

 \subsection{Density of stationary nets and of closed geodesics.}

 We call a stationary net non-trivial if it has at least one vertex of valency $\ge 3$, in other words if it is not a disjoint union of simple closed geodesics.

\begin{qst}[Gromov] {\it Are non-trivial  stationary nets dense on any closed Riemannian surface? In other words, for each non-empty open subset $U$, is there a stationary net intersecting $U$?  } \end{qst}

The answer is evidently positive for every  surface of constant curvature. Indeed, in this case  tangent vectors to closed geodesics are dense in the unit tangent bundle and we can take a finite union of closed  intersecting  geodesics as a stationary net.
 On other manifolds, the answer is less trivial, for example because in
higher dimensions   periodic  geodesics may not intersect.

\section{How to reconstruct a metric from its geodesics.}

A (Riemannian or   semi-Riemannian) metric allows one to construct geodesics. 
Every non-zero vector at every point is tangent to a unique curve from the family of geodesics.
  This section is dedicated to  open problems related to the  following  questions: given a family of curves with this property, are they the geodesics for a metric; and, if there is such  a  metric, is it unique and   how can it be reconstructed? These questions can also be posed for   Finsler metrics.   In this section we assume that the dimension of the manifold is at least two.  
 
  Let us  introduce some notation. By a \emph{path  structure} we will understand  (following \cite{Thomas1925})
  a family of curves $\gamma_{\alpha}(t)$ such that for any point $x$
  and for any vector  $v\in T_xM$, $v\ne 0$,   there exists a unique curve $\gamma$ from this family such that for a certain $t$ we have $\gamma(t)= x$ and $\dot\gamma(t) \in \textrm{span}(v)$.   We  think that the curves are smooth and smoothly depend on the parameters $\alpha= (\alpha_1,\dots,\alpha_{2n-2})$. We may insist that the parametrization of the curves stays fixed, or we may be willing to reparametrize them.
  
  The first step is to  find an affine connection for which the curves are the geodesics. The problem  reduces to a system of linear equations (whether one can actually write down and solve them depends on the form in which the curves are given.)  We outline the main ideas.

        First consider the case in which the parametrization of the curves is fixed.
        If all  the curves $\gamma$ from this path structure   are affinely parameterised geodesics of a connection
         $\nabla=\left(\Gamma^i_{jk}\right)$, then at any point $p\in M$ 
         the Christoffel symbols 
         satisfy the system of equations 
         \begin{equation}\label{christof}
         \ddot \gamma(t)^{i} + \Gamma(\gamma(t))_{jk}^i\dot\gamma(t)^j\dot \gamma(t)^k=0
         \end{equation}
          for all curves  $\gamma$ from the path structure 
            such that $\gamma(t)= p$. This is a system of $n$ 
         linear equations in the $\tfrac{n^2(n+1)}{2}$
          unknowns $\Gamma(p)_{jk}^i$.  Consequently
           $\tfrac{n(n+1)}{2}$   curves $\gamma_\alpha$ from the path structure 
         give us a system of $\tfrac{n^2(n+1)}{2}$  linear  equations in
            $\tfrac{n^2(n+1)}{2}$ unknowns $\Gamma(p)_{jk}^i$. It is an easy exercise  to see that if the velocity vectors of the  curves at the point $p$  are           in    general position,   then this system is  uniquely solvable; by solving it we obtain  
             the Christoffel symbols          at this point. Clearly, these Christoffel symbols should satisfy the equation  
        \eqref{christof} for all curves $\gamma$ from the path structure
         (so generic path structures do not come from a connection).  
          It  depends  on how the curves $\gamma$
         are given   whether it is possibly to check this. For example, if all the curves are given by explicit formulas that depend algebraically on $t$ and on $\alpha$,            then this    is an algorithmically doable but computationally complicated task.

  Note  that the above considerations  show that   reconstruction of a symmetric affine  connection from
  affinely parameterized geodesics is unique. 
     
   Let us now deal with the   reconstruction of a connection from unparameterized curves.  
   Our goal is to find a connection $\nabla=\left(\Gamma^i_{jk}\right)$ for which each curve $\gamma_\alpha$ from our path structure, after an appropriate reparameterization, is
      a geodesic. In this case, a similar idea works. 
   The  analog of \eqref{christof}  is  
    \begin{equation}\label{christof1}
         \ddot \gamma(t)^{i} + \Gamma(\gamma(t))_{jk}^i\dot\gamma(t)^j\dot \gamma(t)^k= c\dot \gamma^i,
         \end{equation}
 where  the unknowns are $\Gamma(p)_{jk}^i$ and  $c$ (though we are interested only in $\Gamma(p)_{jk}^i$). For one curve  $\gamma$  containing $p\in M$ we have therefore $n$ equations in $\tfrac{n^2(n+1)}{2}+ 1$ unknowns $\Gamma(p)_{jk}^i$ and $c$. For two curves $\gamma_1, \gamma_2$ we obtain then $2n$ equations in $\tfrac{n^2(n+1)}{2}+ 2$ unknowns  $\Gamma(p)_{jk}^i, c_1,c_2$ and  so on. We see that for $k>  \tfrac{n^2(n+1)}{2(n-1)}$ 
 curves $\gamma_\alpha$
  from the path structure (passing through the point $p$) we have more equations than unknowns; by
  solving this system (if it is solvable)  we obtain a  connection. See \cite[\S 2.1]{Matveev2012a} for more details.

Note that (as was  already known  to  \cite{Levi-Civita1896} and \cite{Weyl1921}) the solution  $\Gamma_{jk}^i$ of \eqref{christof1}, if it exists,  is not unique:  two connections $\nabla=\left( \Gamma_{jk}^i\right) $ and  $\bar \nabla= \left(\bar \Gamma_{jk}^i\right) $  have the same unparameterized geodesics, 
if and only if   there exists a  $(0,1)$-tensorfield  $\phi_i  $  such    that \begin{equation} \label{c1} 
 \bar \Gamma_{jk}^i  = \Gamma_{jk}^i + \delta_{\ k}^i\phi_{j} + \delta_{\ j}^i\phi_{k}.    
   \end{equation} Thus, the freedom in reconstructing 
   of  a connection is an arbitrary choice of a 1-form $\phi_i$.  

Connections $\nabla = \left(\Gamma^{i}_{jk}\right)$ and  $\bar \nabla = \left(\bar \Gamma^{i}_{jk}\right)$  related by \eqref{c1} are called \emph {projectively equivalent}; projective equivalence of two connections means that they have the same geodesics considered as unparameterized curves.

Let us now touch on the   question of whether/how 
 one can reconstruct a metric (Riemannian or of arbitrary signature)  from a path structure. 
 As we explained above, it is relatively easy to reconstruct an affine  connection (resp.\ 
 a  class of projectively equivalent affine connections)   from affinely (resp. arbitrary) 
 parameterized geodesics, so 
  the actual question  is how to reconstruct a metric from its affine connection (resp.\  a class of projectively equivalent  affine connections), 
  when it is possible, and what is the freedom. Let us first discuss how/whether  it is possible to reconstruct a metric parallel with respect to a given symmetric affine  connection.

  This   question is  well-studied; see for example the answers of Bryant and Thurston  in \cite{Bry-Thu2011}. A  theoretical answer is as follows: the affine connection determines the holonomy group.  The existence of a metric with a given affine connection is equivalent to the  existence of a non-degenerate bilinear form preserved by the holonomy group, see e.g.\ \cite{Schmidt1973}.   
  A practical test for the existence of the metric is as follows: the connection allows to construct the curvature tensor $R^i_{\ jk\ell}$, and if the connection is the Levi-Civita connection of a metric, then this metric satisfies the following equations:
\begin{equation}\label{c2} 
R^s_{\ jk\ell}g_{si}= - R^s_{\ ik\ell}g_{sj}, R^s_{\ jk\ell}g_{si}=   R^s_{\ \ell ij}g_{sk}. 
\end{equation}
These equations are essentially the algebraic symmetries of the Riemann curvature tensor: the first one corresponds to  $R_{ijkm}= - R_{jikm}$, and the second corresponds to  $R_{ijkm}= R_{kmij}$.  
 One should view these equations as linear equations in the unknowns $g_{ij}$; the  number of equations is  bigger than the number of unknowns so it is expected that the system has no non-zero or non-degenerate solution 
  (and in this case there exists no metric compatible with this connection). 
  In many cases  the solution is unique (up to  multiplication by a conformal coefficient) and in this case we already have the conformal class of the metrics. Now,  having the conformal class of the metric we have the  conformal class of the volume form and it is easy to reconstruct the metric using the condition that 
 the volume form is parallel, which   immediately reduces to  
 the condition that a certain 1-form is closed; see also \cite{Mat-Tra2014}.
 
 Note that  if the equations \eqref{c2} do not give enough information one could consider their ``derivatives''  
 $$ 
R^s_{\ jk\ell, m}g_{si}= - R^s_{\ ik\ell,m}g_{sj}, R^s_{\ jk\ell,m}g_{si}=   R^s_{\ \ell ij,m}g_{sk} 
$$
  which gives us  again a huge system of equations in the same unknowns $g_{ij}$. If necessary one then 
   considers higher order derivatives until there is enough information. The general theory says that in the analytic category the existence of a non-degenerate solution of the resulting system of equations  implies the existence of a metric whose Levi-Civita connection is the given one. In the non-analytic setting, however,  there exist $C^{\infty}$ counterexamples.

  Let us now discuss the uniqueness of the reconstruction of a metric from its affinely parameterized geodesics,   modulo the  multiplication of the metric by a constant. A generic  metric can be  reconstructed from its  affinely parameterized geodesics; this follows from the observation that for a generic metric $g$ 
   the system \eqref{c2} has  only the trivial solutions $\textrm{const} \, g$.   So the question is for which metrics the reconstruction is not unique; 
i.e.\ which non-proportional 
   metrics have the same Levi-Civita connection.  
   Locally, the answer is known.    For Riemannian metrics, it was understood already by Cartan and Eisenhart \cite{Eisenhart1923}.
      For metrics of arbitrary signature, the answer is  in the recent papers  \cite{Boubel2015a,Boubel2015b}.
   Both results are a more-or-less  explicit local description  of all 
   metrics having  the same Levi-Civita connection; such metrics are called \emph{affinely equivalent}.    
   
    In the Riemannian case, a global (i.e., when the manifold is closed or complete)
     analog of the result of Cartan and Eisenhart is due to \cite{DeRham1952}, where it was proved   that the universal cover is the direct product of Riemannian  manifolds. 

In the pseudo-Riemannian case the situation is more complicated.   Consider two  affinely equivalent metrics $g$ and $\tilde g$ and the (1,1) tensor $A$ given by  the condition $\tilde g ( A(\cdot), \cdot) = g(\cdot,\cdot)$.  It is self-adjoint with respect to both metrics and 
      parallel with respect to  the Levi-Civita connection.   
       Its generalized eigenspaces  therefore generate mutually orthogonal    totally geodesic foliations the sum of whose dimensions is the dimension of the manifold. Under the additional  assumption that the manifold is geodesically 
      complete one can show, see e.g.\ \cite{Wu1967}, that the universal cover of the manifold splits as the direct product of the (universal covers of the) 
      leaves of these foliations with the restriction of the metrics to them, 
     $$\tilde M  = M_1\times \dots \times M_k,\quad  g= g_{|M_1}+ \cdots + g_{|M_k}, \quad \tilde g= \tilde g_{|M_1}+ \cdots + \tilde g_{|M_k}.$$ 

The most interesting case is when $A$ has only one eigenvalue and since it is clearly a constant one can assume that it is zero.   In this case, all known examples on closed manifolds  are when $g$ admits  light-like parallel vector fields, and $A$ is built with the help of these  vector fields. 
      We therefore suggest, as an interesting  unsolved problem, the following question.   
   
   \begin{qst}
 {\it   Suppose  a    closed manifold $(M,g)$ admits a non-zero 
    $(1,1)$-tensor field that is self-adjoint with respect to $g$,  parallel and nilpotent. Does this manifold or its double cover  admit a non-zero light-like parallel vector field?}
   \end{qst}

   Also in the Finsler case parameterized geodesics determine the connection.  Note though that in the Finsler geometry the connection is non-linear and actually there are more than one popular and widely used connections.  We will call two Finsler   
  metrics \emph{affinely equivalent}  if any   
  geodesic of the first metric   (considered as a curved parametrized  such that the length of the velocity vector is a constant) is a geodesic of the second metric;  this is equivalent to  the Berwald connections of the first and the second metric coinciding in the natural sense.

   \begin{problem}{\it 
Describe all affinely equivalent Finsler  metrics.}   \end{problem}

  A standard example of affinely equivalent non-proportional 
 Finsler metrics can be constructed in the class of Berwaldian metrics. One can also construct easy examples on direct product of Finsler manifolds.  The
 simplest version of the above question is to \emph{construct  non-Berwaldian examples of affinely equivalent non-proportional Finsler metrics which are not  product metrics.}

 The ``unparameterized'' versions of these problems for Riemannian  and pseudo-Riemannian metrics are the subject of the recent survey \cite{Matveev2012a}. Roughly speaking, the situation  is similar to the one in the ``parameterized'' case: locally, a general strategy for reconstructing a metric is understood: the existence of a metric compatible with a path structure   is equivalent to the existence of (non-degenerate) parallel sections of a certain tensor bundle \cite{Eas-Mat2008}. A connection on this tensor bundle is constructed by the projective structure constructed by the path structure.  Thus, 
   a theoretic answer is  to see whether the holonomy group of this connection 
   preserves a certain non-degenerate element of the fiber,  and a practical  method is to construct the curvatures of the connection and look for elements of the fiber compatible with the curvature. As an interesting open
   problem 
  we suggest: 
  
  \begin{problem} \label{prob:br}  {\it 
  Construct a  system of scalar invariants  of a projective structure 
  that vanish if and only if there exists (locally, in a neighborhood of almost every point) a metric compatible with a given projective structure.}   
  \end{problem}
    
   In dimension two, the problem was solved in \cite{Br-Du-Eas2009}, and the system of invariants is quite complicated --- the simplest invariant has degree 5 in derivatives.  It is  possible that in higher  dimensions the system of invariants could be easier in some ways, since in this case the PDE-system  corresponding to the existence of a metric for a projective structure has a higher degree of overdeterminacy. In particular fewer differentiations are needed to construct the first obstruction to the existence of a metric class. See the recent paper \cite{Dun-Eas2014}.

   This  problem  is of course much more difficult in the Finslerian case. 
   In dimension $\ge 3$, not every  path structure is   metrizable, and actually this result is micro-local, see \cite{Douglas1941}. This topic is one of the questions  
    studied  within the topic called {\it inverse problem of calculus of variations}. 
    
 Let us now consider    dimension 2. First observe that if the Finsler metric is irreversible, then its forward and backward geodesics do not necessarily coincide. One  therefore needs to extend the definition of the path structure: by  \emph{irreversible path  structure} we will understand   
  a family of curves $\gamma_{\alpha}(t)$ such that for any point $x$
  and for any vector  $v\in T_xM$, $v\ne 0$,   there exists a unique curve $\gamma$ from this family such that the for a certain $t$ we have $\gamma(t)= x$ and $\dot\gamma(t) $ is proportional to $ v$ with a positive coefficient.    Path structures in the sense of the initial definition  from the beginning of this section will be called  \emph{reversible path  structures}.  
   
  An easy example of an irreversible path structure is provided by circles of radius $1$  on the Euclidean plane  $\mathbb{R}^2$  with  clockwise orientation.  This path structure is locally  Finsler-metrizable, 
  that is there exists a Finsler metric such that  all 
  curves of the   path structure are its geodesics.  
  
  It is known that, in dimension 2,  locally, each reversible  path structure  is metrizable \cite[Theorem B]{Alv-Ber2010}. 
   
   \begin{qst}
 {\it  In dimension $2$, is each irreversible path structure  metrizable? 
   }
   \end{qst} 
   
   We actually ask this question locally. Globally, and on closed manifolds, the answer is negative. A simple example provides  the torus 
    $\mathbb{R}^2/\mathbb{Z}^2$ with  the irreversible path structure from the example above, when the curves from the path structure are circles of radius $1$.  Indeed, in this case not all pairs of points on the universal cover can be connected by a geodesic.

   \weg{The simplest version of the question above is:   
     \begin{qst} 
 {\it  Consider  $\mathbb{R}^2$ with  the Lorentzian flat 
 metric $dx^2 - dy^2$.   Does there exist, locally, a Finsler metric such that each curve of curvature $1$ is a geodesic?}
   \end{qst} 
   
   Recall that  curves of constant curvature equal to $1$ are  the hyperbolas  given by $x^2-y^2=\pm 1$ and their images under isometries of the metric. 
   
   The analog of this question for positive definite metrics of constant curvature has a positive answer since curves of constant curvature are actually magnetic geodesics with $\omega$ being  the volume form of $g$, and therefore are geodesics of a certain Randers metric, which was calculated explicitly in 
    \cite{Cra-Mes2014}.  }

Let us now discuss the freedom in reconstructing the (Riemannian or pseudo-Riemannian)   metric from its unparameterized geodesics; an equivalent question is how many different metrics can have the same geodesics considered as unparameterized curves.  Of course, the metrics $g$ and $\textrm{const} \cdot g$ have the same geodesics. In \cite{Matveev2012a} it was shown that, for a $k$ depending on the dimension, a $C^k$-generic metric has the property that any projectively equivalent metric is proportional to it with a constant coefficient.   In fact, as  follows from \cite{Kru-Mat2016}, in all dimensions   $k=2$ is enough.  
 
There exist local (the first examples are due already to Lagrange, Beltrami and Dini) 
and global examples of non-proportional projectively equivalent  metrics. Locally, in the Riemannian case, a complete description of projectively equivalent metrics  is due to Levi-Civita \cite{Levi-Civita1896}    and in arbitrary signature is due to \cite{Bol-Mat2015}. Globally, in the Riemannian case, the situation is also pretty clear,  see e.g.\ \cite{Mat-Top1998,Matveev2003a,Matveev2003b},  but if the metrics have arbitrary signature, virtually nothing is known. A general problem is to understand the topology of closed manifolds admitting non-proportional projectively equivalent metrics of indefinite signature   and  as the simplest  version of this problem we suggest

   \begin{qst} {\it Can a 3-dimensional sphere admit two non-proportional  projectively equivalent
   metrics of Lorentzian signature?}  
   \end{qst}

In the Finslerian case, the freedom of reconstructing metric from its geodesic can  be much higher. For example, adding 
any sufficiently small closed  1-form to the  Finsler metric 
 does not change  geodesics of this metric. This transformation of the metric was called  a \emph{ trivial projective change} in \cite{Matveev2013}.

  \begin{qst} {\it In dimension $\ge 3$,  for a generic Finsler metric $F$, is any  metric projectively equivalent to $F$   
    locally proportional to a trivial projective change of the metric? }
   \end{qst} 

In dimension $2$, generic Finsler metrics admit projectively equivalent metrics which do not come from a trivial projective change.

\subsection*{Acknowledgements.} We thank  all the participants, 
especially Victor Bangert, Charles Boubel, Nancy Hingston, Gerhard Knieper, Yiming Long, Gabriel Paternain, Mark Pollicott, Regina Rotman, St\'ephane Sabourau, Benjamin Schmidt, for active participation in the problem sessions. We  also thank non-participants who have helped: Alberto  Abbondandolo,  Juan-Carlos Alvarez Paiva, Ivan Babenko,   Florent 
Balacheff, Helga Baum, Michael Bialy, Dima Burago,  Maciej Dunajski, Misha Gromov, Larry Guth, Sergei Ivanov,  Bruce Kleiner, Alan Reid, Stefan Suhr.

Most of the writing of the paper was done during four  visits by KB to the  Friedrich-Schiller-Universit\"at Jena, supported by DFG (SPP 1154 and GK 1523).


\end{document}